\documentclass[a4paper,11pt]{article}
\usepackage[a4paper]{geometry}
\usepackage[utf8]{inputenc}

\title{Whitehead algorithm for generalized Baumslag-Solitar groups}
\author{Chloé Papin}

\usepackage{graphicx}
\usepackage{amssymb}
\usepackage{tikz}
\usetikzlibrary{shapes}

\usepackage{amsmath}
\usepackage{amsthm}
\newtheorem{theo}{Theorem}[section]
\newtheorem{theointro}{Theorem}
\newtheorem{prop}[theo]{Proposition}
\newtheorem{propintro}[theointro]{Proposition}

\newtheorem{lem}[theo]{Lemma}
\newtheorem{lemintro}[theointro]{Lemma}
\newtheorem{cor}[theo]{Corollary}
\newtheorem{corintro}[theointro]{Corollary}

\theoremstyle{definition}
\newtheorem{defi}[theo]{Definition}
\theoremstyle{remark}
\newtheorem{rema}[theo]{Remark}
\newtheorem{remas}[theo]{Remarks}
\newtheorem{ex}[theo]{Examples}

\def \N{\mathbb N}
\def \Z{\mathbb Z}

\def \D{\mathcal D}
\def \A{\mathcal A}
\def \C{\mathcal C}
\def \E{\mathcal E}
\def \G{\mathcal G}

\DeclareMathOperator{\Out}{Out}

\DeclareMathOperator{\axe}{Axis}
\DeclareMathOperator{\stab}{Stab}

\DeclareMathOperator{\id}{id}
\DeclareMathOperator{\Amin}{\mathcal{A}_{\text{min}}}

\begin{document}

\maketitle
\newcommand{\localhost}{.}
 
 \begin{abstract}
 In analogy with the free factors of a free group we define \emph{special factors} of Generalized Baumslag-Solitar (GBS) groups as non-cyclic subgroups which appear in splittings over infinite cyclic groups. We give an algorithm which, given a GBS group $G$ and an element $g \in G$, decides whether there exists a special factor $H < G$ such that $g \in H$. This algorithm is analogous to an algorithm by Whitehead for free groups. Furthermore we prove that given $g \in G$ there exists a unique minimal special factor containing $g$ and give an algorithm which finds it.
 \end{abstract}
 
 \bigskip

Baumslag-Solitar groups, defined by $BS(p,q)=\langle a, t | ta^pt^{-1}=a^q \rangle$, are a family of one-relator groups which were introduced to give examples of non-hopfian groups in \cite{BaumslagSolitarSomeTwoGenerator}. Generalized Baumslag-Solitar (GBS) groups are a wider family composed of all fundamental groups of finite graphs of infinite cyclic groups. The isomorphism problem between GBS groups (i.e. determining if two graphs of groups define the same GBS group) is not known apart from a few special cases. Clay and Forester showed in \cite{ClayForesterOnTheIsomorphism} that it is solvable when first Betti number is at most one. Levitt (\cite{LevittAutomorphisms}) defined a particular class of GBS groups with nice outer automorphism group for which the isomorphism problem is also solvable.

In fact, the automorphism group of a GBS group $G$ can be very complicated in general and depends a lot on $G$. For example, $Out(BS(2,3))$ is finite but $Out(BS(p,pn))$ for $p,n \in \Z$ is not finitely generated. 

This was shown by Collins and Levin in \cite{CollinsLevinAutomorphisms} with algebraic methods and later by Clay in \cite{ClayDeformation} using geometric methods.

By Bass-Serre theory GBS groups have a natural action on trees which is a motivation for using geometric methods. The study of automorphisms of free groups and outer space give some inspiration. Given a GBS group $G$, there is an analogue for outer space, which is called \emph{deformation space} (\cite{ForesterSplittings}). It can be defined as the set of all minimal actions of $G$ on trees with infinite cyclic stabilizers, up to $G$-equivariant isometry 
and $\Out(G)$ acts on it by precomposition of the action. 

An automorphism of the free group $F_N$ is called \emph{fully irreducible} if it does not have any periodic conjugacy class of free factor of $F_N$. Such an automorphism $\phi$ admits a \emph{train track} representative: there exists a free minimal action of $F_N$ on a simplicial metric tree, which defines a translation axis for the action of $\phi$ on outer space. Train tracks representatives are a powerful tool to study automorphisms of $F_N$.

In some sense one can define fully irreducible automorphisms of a GBS group $G$: they are the automorphisms whose powers do not preserve any conjugacy class of \emph{special factors}, the analogue of free factors in this context. 
Margot Bouette showed in \cite{BouetteThese} that all automorphisms of $BS(p,pn)$ are reducible, which is surprising: 
one would have suspected that a generic automorphism would not preserve any special factor. However, she proved that there is a different deformation space invariant under $Out(G)$ on which irreducible automorphisms exist, always admit train tracks and act on the modified deformation space with positive translation length. In more general GBS groups it is not known if train tracks always exist for fully irreducible automorphisms. Under some restrictions on the GBS group though (ex: no nontrivial integer modulus, \cite{ForesterSplittings}), we can assume that the dimension of the deformation space is finite. Then \cite[Theorem 50]{MeinertTheLipschitzMetric} applies and proves the existence of train tracks.

This sort of problem is a motivation to understand \emph{special factors}. 
Here is a way to define them. Identifying the free group $F_N$ with the fundamental group of some graph $\Gamma$, the fundamental group of any subgraph of $\Gamma$ is a free factor of $F_N$. Moreover, allowing to vary the graph and identification, any free factor of $F_N$ can be obtained this way. Similarly, a special factor of a GBS group $G$ is the fundamental group of some subgraph of a graph of groups $\Gamma$ where $\Gamma$ is a graph of cyclic groups with fundamental group $G$. To avoid degenerate cases, we do not consider cyclic subgroups as special factors. A noticeable difference is the fact that many different graphs of groups may appear and give distinct special factors, while all free factors of free groups may be seen in roses. Therefore, while free factors of same rank in a free group $F_N$ are isomorphic and in the same orbit under $Aut(F_N)$, special factors of a GBS group may be a lot more diverse and there may be infinitely many orbits of special factors under $Aut(G)$. This is the case for $G\simeq BS(2,4)$.

The Whitehead algorithm in the free group solves the following problem: given $g,h \in F_N$, is there $\phi \in Aut(F_N)$ such that $\phi(g)=h$ ? A weaker form of this algorithm can decide, given $g \in F_N$, if $g$ is simple, i.e. whether there exists a proper free factor containing $g$. 
The first problem seems difficult for GBS groups since the automorphism group is a lot more complicated and may not be finitely generated. 
Here we solve the analogue of the second problem in a GBS group $G$ by giving an algorithm deciding if an element $g \in G$ is contained in a special factor. We also prove that there exists a unique minimal special factor containing $g$ and give a further algorithm which finds this minimal factor.

 \bigskip
 
Before stating the results, let us give some useful background.
Let $G$ be a GBS group. It admits an action on a locally finite tree $T$ with cyclic edge and vertex stabilizers. 
Unless $G$ is isomorphic to $\Z, \Z^2$ or the fundamental group of a Klein bottle, all trees with cyclic stabilizers have the same elliptic subgroups.
Let $\D$ be the set of $G$-trees with cyclic edge and vertex stabilizers. We call it the \emph{cyclic deformation space}. 

In the main part of the paper we will consider \emph{restricted deformation spaces}, where we impose an extra condition on the edge groups of the trees. In this introduction we present the result for $\D$ the cyclic deformation space for simplicity but we actually prove a slightly stronger version.

A \emph{loxodromic} element in a tree $T$ in $\D$ is an element acting on a tree with no fixed point. This property does not depend on the choice of $T \in \D$, so we may refer to loxodromic elements of the group $G$. A loxodromic element of $G$ is \emph{simple} if there exists a proper special factor containing it. In order to understand whether an element is simple, we use \emph{Whitehead graphs}. The Whitehead graph $Wh_T(g,v)$, where $v$ is a vertex, has vertices corresponding to the edges of $T$ with origin $v$. Two vertices are linked by an edge whenever a translate of the axis of $g$ in $T$ takes the turn between the two edges corresponding to the vertices. We say that the graph has an \emph{admissible cut} when it is disconnected or has a cut vertex.

Whitehead's lemma was originally published in \cite{WhiteheadOnEquivalentSets} for free groups. We have the following result, adapted from a version of Whitehead's lemma from \cite{GuirardelHorbezAlgebraic}.

\begin{theointro} \label{theointro:whitehead}
 Let $g \in G$ a loxodromic element. If $g$ is simple then for all $T \in \D$ there exists a vertex $v \in T$ such that $Wh_T(g,v)$ has an admissible cut.
\end{theointro}

Let us take an element $g \in G$ and a tree $T\in \D$. We consider $T$ as a marked graph of groups $\Gamma = T/G$. In case $Wh_T(g,v)$ has an admissible cut, either $g$ is represented in the fundamental group of some proper subgraph  of $\Gamma$, so it is simple, or we can perform some transformation of $\Gamma$ in order to obtain a new tree. We can apply the theorem above to the result of the transformation. The point is that this process will eventually stop, so we deduce an algorithm to determine whether $g$ is simple:

\begin{corintro} \label{corintro:algo}
 There is an algorithm taking as input
 \begin{itemize}
  \item a GBS group $G$ given as a graph of groups
  \item a loxodromic element $g\in G$
 \end{itemize}
 which decides whether $g$ is simple, and if it is,  returns a proper special factor containing $g$.
\end{corintro}

The point of the algorithm will be to transform the graph until one of its proper subgraphs contains the axis of $g$. Note that Whitehead's version for free groups (see \cite[Chapter I, section 4]{LyndonSchuppCombinatorial}) has the same goal but uses specific automorphisms to transform a rose. Here we prefer to use the approach given in \cite{GuirardelHorbezAlgebraic} instead. It consists in \emph{unfolding} the tree associated to the graph, that is to say, perform the inverse of folds. While in Whitehead's approach the word length decreases and the graph keeps the same volume, here the word length remains constant while the volume of the graph increases, eventually reaching a point where $g$ avoids some edge in the graph.

The main reason for which we choose this approach is that there is no preferred graph of groups in general, and the graph of group showing the special factor containing $g$ might be very different from the initial graph of groups used to define $G$. There might be no automorphism between the factor containing $g$ and a subgraph of the initial graph.

We prove a stronger version of Theorem \ref{theointro:whitehead} and Corollary \ref{corintro:algo} which applies to a finite collection of elements of $G$ rather than a single element. Actually the version for collections decides whether the elements of a collection belong to a \emph{system of special factors}, which are special factors which are either conjugate or disjoint in some sense.

\bigskip

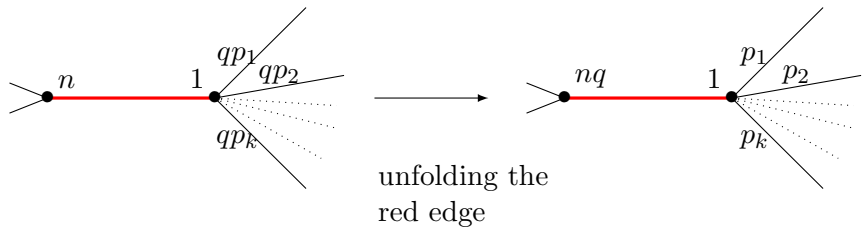
\begin{figure}
 \centering
 \begin{tikzpicture}[]
%\draw [dotted] (0,-6) rectangle (11.5,6);

%le haut
\begin{scope}[]

	\coordinate (w) at (0.7,0);
	\coordinate (v) at (2.9,0);

	%arete milieu
    \draw [very thick, red] (w)-- (v);
	\draw  (w) node {$\bullet$} node [above right] {$n$};
    \draw (v) node {$\bullet$} node [above left] {$1$};

	%arêtes droite
	\draw (v)  -- +(1.2,1.2)  node [near start, above] {$qp_1$};
	\draw (v)  -- +(1.7,0.3)  node [midway, above=-0.1cm] {$qp_2$};
	\draw (v)  -- +(1.2,-1.2)  node [near start, below] {$qp_k$};
	\draw[dotted] (v)  -- +(1.6,-0.1) ;
	\draw[dotted] (v)  -- +(1.6,-0.4); 
	\draw[dotted] (v)  -- +(1.4,-0.8); 
	%arêtes gauche
	\draw (w) -- +(-0.5,0.2);
	\draw (w) -- +(-0.5,-0.2);

\end{scope}
%la flèche
	\draw[-latex] (5,0) -- (6.5,0);
	\node [text width=2.3cm] (leg) at (6.2, -1.3)  {unfolding the red edge};
% dessin du bas 
\begin{scope}[xshift=6.8cm]

	\coordinate (w) at (0.7,0);
	\coordinate (v) at (2.9,0);

	%arete milieu
    \draw [very thick, red] (w)-- (v);
	\draw  (w) node {$\bullet$} node [above right] {$nq$};
    \draw (v) node {$\bullet$} node [above left] {$1$};

	%arêtes droite
	\draw (v)  -- +(1.2,1.2)  node [near start, above] {$p_1$};
	\draw (v)  -- +(1.7,0.3)  node [midway, above=-0.1cm] {$p_2$};
	\draw (v)  -- +(1.2,-1.2)  node [near start, below] {$p_k$};
	\draw[dotted] (v)  -- +(1.6,-0.1) ;
	\draw[dotted] (v)  -- +(1.6,-0.4); 
	\draw[dotted] (v)  -- +(1.4,-0.8); 
	%arêtes gauche
	\draw (w) -- +(-0.5,0.2);
	\draw (w) -- +(-0.5,-0.2);

\end{scope}

\end{tikzpicture}
 \caption{Unfolding of type II (i.e. the unfolded edges belong to the same orbit) as seen in the graph of groups} \label{fig:lesPliagesDeLintro}
\end{figure}

In the process of showing that the algorithm stops, we need to show that the volume of the graph increases. It does not increase at each step and in fact, if the transformation performed is an \emph{unfolding} of edges in the same orbit (see Figure \ref{fig:lesPliagesDeLintro}), the volume remains constant. The following result prevents the graph from being transformed indefinitely without increasing its volume.

A \emph{sequence of unfoldings of type II} is a finite or infinite sequence $T_n \to \dots \to T_1 \to T_0$ where the $T_i$ are trees in $\D$ and $T_{i+1} \to T_i$ is a type II fold, that is, a fold of edges in the same orbit, with origin and endpoint in different orbits. Then no such sequence can be infinite:

\begin{lemintro}
 Suppose that $G$ is not a solvable Baumslag-Solitar group.  Any sequence of unfoldings of type II is finite. 
\end{lemintro}

 Whitehead graphs, which are useful for the algorithm, can be computed algorithmically. In \cite{BeekerMultiple} methods for algorithmic computations in GBS groups are given. In section \ref{sec:algo} we explain the algorithm in detail.

\bigskip

\begin{figure}
 \centering
 \begin{tikzpicture}[]
 
  \draw (1,0) node {$\bullet$} arc (0:360:1) node [above left] {$2$} node [below left] {$1$} node [midway, sloped] {$>$} node [midway, left] {$t$};
  \draw (1,0) node [above right] {$2^k$}  node [below right] {$\langle b \rangle $} -- (3,0) node {$\bullet$} node [above left] {$2$} node [right] {$\langle a \rangle$};

    \draw (8,0) node [text width=8cm]  {Graph of groups for $BS(2,4)$};

    \draw[-latex]  (2,-0.7) -- (2,-2) node [right, midway] {Expansions};

  \begin{scope}[yshift=-3.5cm, xshift=0]
      \draw (1,0) node {$\bullet$} arc (0:360:1) node [above left] {$2$} node [below left] {$1$} node [midway, sloped] {$>$} node [midway, left] {$t$};
      \draw (1,0) node [above right] {$2$}  -- (2.5,0) node {$\bullet$} node [above left] {$1$} node [above right] {$2$}  -- (4,0) node {$\bullet$} node [above left] {$1$} node [above right] {$2$}   -- (8.5,0) node [midway, above] {$\dots$} node {$\bullet$} node [above left] {$1$} node [above right] {$2$} -- (10,0) node {$\bullet$} node [above left] {$2$};

      \draw (1,-1) |- (5.6,-1.3) node [below] {$k$ edges} -| (10,-1);

    \draw[dashed, blue] (8, 0.7) rectangle (10.3,-0.4);
    \draw[dashed, green] (3.5, 0.9) rectangle (10.5,-0.6);
    \draw[dashed, red] (2,1.1) rectangle (10.7,-0.8);
  \end{scope}

\end{tikzpicture}
 \caption{A construction of arbitrary long special factor sequences in $BS(2,4)$. The upper graph represents $BS(2,4)$ with the presentation $\langle a,b,t| tbt^{-1}=b^2, a^2=b^{2^k} \rangle$. Expansions from this graph lead to the graph below, which has many edges. The rectangles show some special factors which can be read in the graph below. Here we can construct a sequence of $k$ nested special factors.} \label{fig:chaineDecroissante}
\end{figure}
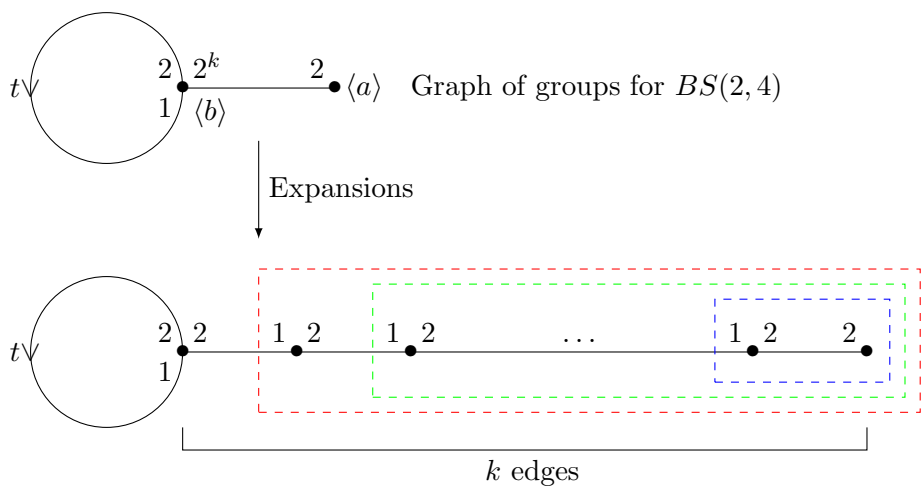

As a special factor is itself a GBS group, we may want to iterate the algorithm in the following way. Given $g \in G$ belonging to a special factor $H \in G$, we can apply the algorithm to $H$ in order to find a smaller special factor, again and again. In the case of free groups, the rank of the free factor decreases, so eventually we find a free factor which does not have any smaller free factor containing the element. In the case of GBS groups the rank does not always decrease, but we do have a complexity $\C$ which decreases strictly when passing to a proper special factor. Note that in some GBS groups including $BS(2,4)$ there exist arbitrary long chains of decreasing special factors (see Figure \ref{fig:chaineDecroissante}). For the definition of $\C$, see section \ref{sec:suitesDecroissantes}.

\begin{propintro}
 Let $G$ be a GBS group and $H$ a special factor of $G$. Then $\mathcal{C}(H) < \mathcal{C}(G)$.
\end{propintro}

From this we deduce the existence of a minimal special factor containing $g$, which is in fact unique:

\begin{theointro}\label{theointro:factMin}
 The set of special factors relative to $\D$ which contain a given loxodromic element $g$ admits a smallest element for inclusion.
\end{theointro}

Using the previous algorithm repeatedly we get the following:
\begin{theointro} \label{theointro:algo}
 There exists an algorithm taking as input a GBS group $G$ as a marked graph of groups and a hyperbolic element $g$, and which outputs the minimal factor containing $g$ as a subgraph of a marked graph of groups for $G$.
\end{theointro}

Just like Theorem \ref{theointro:whitehead}, Theorem \ref{theointro:algo} also applies to finite collections of elements of $G$. In that case it outputs the minimal system of special factors containing the collection.

\bigskip

The paper is structured as follows. In section \ref{sec:definitions} we introduce basic notions about GBS groups and prove that no sequence of unfoldings of type II can be infinite. In section \ref{sec:facteursSpeciaux} we prove Theorem \ref{theointro:whitehead}. Section \ref{sec:algo} proves Theorem \ref{theointro:algo} : it gives of an algorithm to determine whether there exists a special factor containing a particular group element. Most of it consists in proving why all steps needed to check that an element is simple are algorithmic. Since it is more technical, the proof of the algorithm may be skipped at first reading. Finally we prove theorem \ref{theointro:factMin} in section \ref{sec:suitesDecroissantes}.

\section{Definitions} \label{sec:definitions}

\subsection{Trees, elliptic groups, deformation spaces}

 We refer to \cite{SerreArbresAmalgames} for basic notions on graphs of groups.
 A \emph{graph} is a set of vertices, along with a set of oriented edges. An edge comes with applications $o$ and $t$ which associate its initial and terminal vertex to an edge, respectively. There is a fixed-point-free involution $e \mapsto \bar e$ with $t(\bar e ) = o(e)$.

 A \emph{graph of groups} is a graph $\Gamma$ such that every vertex $v$ is labelled with a group $G_v$, every edge $e$ is labelled with a group $G_e = G_{\bar e }$ and for every oriented edge $e$ of $\Gamma$, there is a given monomorphism $\phi_e : G_e \to G_{t(e)}$.
 
\begin{defi}
  A \emph{generalized Baumslag-Solitar group} (\emph{GBS} group) is the fundamental group of some finite graph of groups of which each edge or vertex group is infinite cyclic.
  
  A GBS group is \emph{non-elementary} if it is not isomorphic to one of the following: $\Z$, $\Z^2$ or the fundamental group of the Klein bottle $\langle a, t | tat^{-1}=a^{-1} \rangle = \langle a,b | a^2=b^2 \rangle$.
\end{defi}

\begin{remas}
 \begin{itemize}
  \item Generally the graph of (cyclic) groups for a GBS group is not unique, in fact there may be infinitely many such graphs.
  \item From a graph of groups one can deduce a presentation for the group. Since vertex groups are cyclic, GBS groups are finitely presented. See \cite{LevittAutomorphisms} for some detail on GBS groups.
 \end{itemize}
\end{remas}

\begin{defi}
 A \emph{labelled graph} is a graph of which each oriented edge $e$ carries a label $\lambda (e) \in \Z \setminus \{0\}$ at its origin.
\end{defi}

In graphs of cyclic groups, choosing a generator $a_v$ for each vertex group and $a_e$ for each edge group gives an identification with $\Z$ and each inclusion $G_e \to G_v$ is given by a multiplication by a nonzero integer. Therefore any such graph of groups can be described by a labelled graph.

\begin{defi} 
 An edge in a graph is a \emph{loop} if its endpoints are equal. 
\end{defi}

\bigskip

Let $G$ be a non-elementary GBS group.

\begin{defi}
 A \emph{$G$-tree} is a simplicial tree endowed with a minimal action of $G$ by simplicial isomorphisms without inversion of edges.
 
 We endow all trees with the combinatorial metric, that is, all edges have length $1$.
\end{defi}

\begin{rema}
 Let $T$ be a $G$-tree and let $g$ be a loxodromic element in $T$, that is, an element acting on $T$ with no fixed point. The \emph{translation length} of $g$ is 
 \[
  \|g\|_T := \min_{x \in T} d(x,g x)
 \]
 and it is equal to of edges in a fundamental domain of the axis of $g$ in $T$.
\end{rema}

Let $T$ be a $G$-tree. Let $e$ be an edge of $T$. We denote the stabilizer of $e$ by $G_e$ (and this subgroup fixes both endpoints of $e$ since the action is without inversion). Similarly we denote by $G_v$ the stabilizer of a vertex $v$.

\begin{defi} 
 Let $T$ be a $G$-tree and $H$ a subgroup of $G$. The subgroup $H$ is \emph{elliptic} (resp. \emph{bi-elliptic}) in $T$ if it fixes a vertex (resp. an edge).
\end{defi}

\begin{defi}
 Two trees $T,T'$ are in the same \emph{deformation space} if they share the same elliptic subgroups. If a subgroup $H < G$ fixes a point in $T$, it must then fix a point in $T'$, and conversely. 
\end{defi}

\begin{rema} \label{rem:domine}
 Equivalently, $T$ and $T'$ belong to the same deformation space if there exists $G$-equivariant applications $T \to T'$ and $T' \to T$ (see \cite{GuirardelLevitt07} for detail).
 
 Apart from $\Z, \Z^2$ and the fundamental group of the Klein bottle, all actions on trees dual to descriptions of $G$ as graphs of cyclic groups belong to the same deformation space (\cite{ForesterSplittings}).
\end{rema}

\begin{defi}
 For a non-elementary GBS group $G$, the \emph{cyclic deformation space} $\D$ is the space of all $G$-trees with infinite cyclic vertex and edge stabilizers.
\end{defi}

 We define a restricted deformation space of the cyclic deformation space, which is smaller than $\D$, like in \cite{GuirardelLevitt07}:

\begin{defi}
 Let $\mathcal A$ be a family of subgroups of $G$, stable by conjugation and by taking subgroups. The \emph{restricted deformation space} $\D^{\mathcal A}$ is the set of $G$-trees in $\D$ whose bi-elliptic subgroups all belong to $\mathcal{A}$.
 
 We refer to subgroups which belong to $\mathcal A$ as \emph{allowed edge groups} in $\D^{\mathcal A}$.
\end{defi}

\begin{rema}
 For a subgroup of $G$, being elliptic is a property which depends only on the deformation space of the tree studied. However in general this does not hold for bi-elliptic subgroups: two trees in the same deformation space may have different bi-elliptic subgroups. This applies to $\D$ and $\D^{\mathcal A}$. 
\end{rema}

\begin{defi}
 Let be $T$ a $G$-tree. Let $F$ be a $G$-invariant subforest of $T$. We define the equivalence relation $\sim_F$ on $T$ as the smallest $G$-invariant equivalence relation such that $x \sim_F y$ whenever there exists some connected component of $F$ containing both $x$ and $y$.
 
 The quotient $T \to T/\sim_F$ is a tree, called the \emph{collapse} of $F$ in $T$.

 The $G$-equivariant application $T \to T/\sim_F$ is called the  \emph{collapse map}.
  
  When $F$ is the $G$-equivariant subforest spanned by an edge $e$, we speak of the \emph{collapse of $e$} instead of the collapse of $F$, and denote it by $T \to T/\sim_e$ 
\end{defi}

\begin{defi}
 Let $T$ be a $G$-tree. An edge of $T$ is called \emph{collapsible} if $T/\sim_e$  is a tree in the same deformation space as $T$. An edge $e$ with endpoints $u$ and $v$ is collapsible if $u$ and $v$ are in different orbits and either $G_e = G_u$ or $G_e= G_v$.
\end{defi}

\begin{rema}
 If there is a $G$-equivariant application  $f: T \to S$ where $T,S \in \D$ such that the image of $e$ is a single point in $S$, then $e$ is collapsible. Indeed $f$ factorizes through $T/ \sim _{e}$. Thus there exists a $G$-equivariant application $T/ \sim _{e} \to T$ (and vice-versa) so $T/ \sim _{e}$ is in $\D$.
\end{rema}

\begin{defi}
 A $G$-tree $T$ is \emph{reduced} if no edge in $T$ is collapsible.
\end{defi}

\begin{remas}
 \begin{enumerate}
  \item This notion can be expressed of in terms graphs of groups, which are more convenient for a computational use. Denote by $\Gamma$ the graph of groups associated to the quotient $T/G$, then $T$ is reduced if whenever an edge morphism in $\Gamma$ into a vertex group is surjective then the edge is a loop.
  \item Reduced trees in a deformation space share the same bi-elliptic subgroups. Any bi-elliptic subgroup in some reduced tree $T \in \D$ is also bi-elliptic in any tree in $\D$, reduced or not. 
 \end{enumerate}
\end{remas}

\begin{defi}
 We call $\mathcal A _{\text{min}}$ the family of subgroups which are bi-elliptic in reduced trees. The \emph{reduced deformation space} $\D_{\text{red}}$ is the restricted deformation space $\D^{\mathcal A _{\text{min}}}$.
\end{defi}

 \begin{defi}
 We call an elliptic subgroup in some tree $T \in \D$ \emph{big} with respect to $\mathcal A$ if it is not in $\mathcal{A}$. 
 \end{defi}

\begin{rema}\label{rem:grosSommets}
 In order to check whether an elliptic subgroup is big for $\Amin$, it suffices to check that it fixes no edge in some arbitrary tree in $\D$. 
 
 In a graph of cyclic groups viewed as a labelled graph, vertex groups such that no label at the vertex is $\pm 1$ are maximal big groups. 
  If the graph of groups is reduced, then maximal big groups are exactly vertex stabilizers with all labels different from $\pm 1$.
\end{rema}

The number of conjugacy classes of maximal big groups is finite (bounded by the number of vertex orbits of some tree in $\D$) and depends only on $G$ and $\D$ (see \cite{GuirardelLevitt07}).

\subsection{Folds, expansions} \label{subsec:folds}

In this section we define folds and expansions, and give a construction for a certain type of expansion.

\begin{defi}
 A $G$-tree $S$ is a \emph{refinement} of  another $G$-tree $T$ if there exists a collapse $\pi : S \to T$, i.e. $T$ is equivariantly isomorphic do $S/\sim_F$ for some $G$-invariant forest $F \subset S$.
 
 We say that $S$ is an \emph{expansion} if additionnally $S$ and $T$ belong to the same deformation space.
\end{defi}

\begin{lem}[Construction of an expansion $T^{H,S}$] \label{lem:constructionEclatement}
 Let $T$ be a $G$-tree and $v$ a vertex of $T$. Let $\E_v$ be the set of edges with origin $v$. Let $H$ be a subgroup of $G_v$ and $S$ a non-empty proper subset of $\E_v$ which satisfy:
  \begin{itemize}
  \item $HS=S$
  \item $\forall g \in G_v \setminus H , \: gS \cap S = \varnothing$.
 \end{itemize} 
 Consider the partition $\E_v = \bigsqcup _{g \in G_v /H} gS \sqcup E'$ where $E'= \E_v \setminus G_v \cdot S$. 
 
 We construct $T^{H,S}$ from $T$ as follows. First we replace $v$ by the star on $G_v/H$ and for each $g \in G_v$ and each edge $e \in S$, we attach $g \cdot e$ to $gH$. Then we attach $E'$ on the centre of the star. Finally we extend this by equivariance to all translates of $v$.
 
 The obtained tree $T^{H,S}$ is minimal, is a refinement of $T$ and belongs to the same deformation space (without restriction on edge groups). Moreover, if $H \in \mathcal A$ and $T \in \D^{\mathcal A}$ then $T^{H,S} \in \D^{\mathcal A}$.
\end{lem}

\begin{proof}
 Let us prove that $T^{H,S}$ is minimal.
 The tree $T^{H,S}$ has no valence 1 vertex: since $S \subsetneq \E_v$ and $S \neq \varnothing$ then either $G_v/H$ has two or more elements, or $E'$ is non-empty. In both cases no vertex has valence 1.
 
 If $T^{H,S}$ were not minimal, it would have a valence 1 vertex: suppose $T^{H,S}$ is not minimal and call $T^{H,S}_{\text{min}}$ the minimal invariant subtree. By cocompacity of the action there exists $v \in T$ such that the distance $d(v, T_{\text{min}})$ is maximal and positive. Then $v$ has valence 1, which gives the contradiction needed.

 Note that collapsing the stars involved in the construction of $T^{H,S}$ allows to recover the original tree $T$, so $T^{H,S}$ is a refinement of $T$. This implies that all elliptic subgroups of $T^{H,S}$ are elliptic in $T$. Conversely, let $G_0$ be a subgroup fixing a vertex $w$ in $T$. If $w \notin G \cdot v$ then its pre-image by the collapse map is a single point so it must be fixed by $G_0$. If $w=v$, its pre-image is the closure of the star. The star is invariant by $G_v$ so its center is a fixed point for $G_v$ which contains $G_0$. This proves  that $T$ and $T^{H,S}$ have the same elliptic subgroups.

 Thus $T^{H,S}$ is in the same deformation space as $T$. Note that the second condition implies that the stabilizer of any edge of the star is conjugate to $H$, hence the last statement.
\end{proof}

To define folds, we rely on
\cite{BestvinaFeighnBounding}. Folds are a classical operation on trees and were first defined by Stallings (\cite{StallingsFoldings}).
The idea behind a fold is to identify two edges of a $G$-tree in an equivariant way to create a new tree.
\begin{defi}[Fold] 
 Let $T$ be a $G$-tree. Let $V$ be a vertex of $T$ and $e_1,e_2$ two edges with origin $v$ such that $e_2 \notin G\cdot \bar e_1$. We define the fold of $e_1$ together with $e_2$ as follows. 
 
 We define the equivalence relation $\sim$ on $T$ as the smallest $G$-invariant equivalence relation satisfying $e_1 \sim e_2$ and $t(e_1) \sim t(e_2)$.

 The result of the fold of $e_1$ with $e_2$ is $T/\sim$ and it is a tree.
\end{defi}

\begin{figure}[h]

\centering
\begin{tikzpicture}
\begin{scope}
\draw (0,0) node [left] {$A$} node {$\bullet$} -- (2,1) node [above] {$B_1$}node {$\bullet$} node [midway,above]{$E_1$};
\draw (0,0)  -- (2,-1) node [above] {$B_2$}node {$\bullet$}node [midway,above]{$E_2$};

\draw[-latex] (3,0)--(4,0);
\draw (5,0) node [left] {$A$} node {$\bullet$} -- (7.5,0) node [above] {$\langle B_1, B_2 \rangle$}node {$\bullet$} node [midway, below]{$\langle E_1, E_2 \rangle$};

\draw (-2,0) node {Type I};
\end{scope}

\begin{scope}[yshift=-3cm]
\draw (0,0) node [left] {$A$} node {$\bullet$} -- (2.5,0) node [above] {$B$}node {$\bullet$} node [midway, below]{$E$};

\draw[-latex] (3,0)--(4,0);
\draw (5,0) node [left] {$A$} node {$\bullet$} -- (7.5,0) node [above] {$\langle B,g \rangle$}node {$\bullet$} node [midway, below]{$\langle E,g \rangle$};

\draw (-2,0) node {Type II};
\end{scope}

\begin{scope}[yshift=-6cm]
\draw (0,0) node [left] {$A$} node {$\bullet$} .. controls (0.7,0.8) and (1.3,0.8) .. (2,0) node {$\bullet$} node [right] {$B$} node [midway, above] {$E_1$};
\draw (0,0)  .. controls (0.7,-0.8) and (1.3,-0.8) .. (2,0)  node [midway, above] {$E_2$};

\draw[-latex] (3,0)--(4,0);

\draw (5,0) node [left] {$A$} node {$\bullet$} -- (7.5,0) node [above] {$\langle B,g \rangle$}node {$\bullet$} node [midway, below]{$\langle E_1,E_2 \rangle$};

\draw (-2,0) node {Type III};
\end{scope}
\end{tikzpicture}
\caption{The three basic types of folds} \label{fig:typesPliages}
\end{figure}
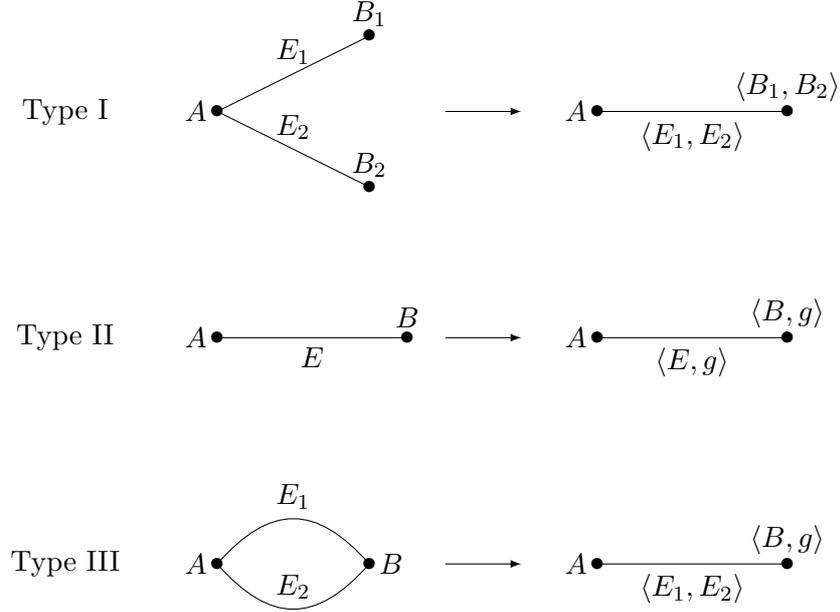

Up to subdividing some edges, any fold boils down to a sequence of folds of the three types illustrated by Figure \ref{fig:typesPliages} (see \cite{BestvinaFeighnBounding} for more details on fold types), which presents folds as seen in the quotient graph, depending on whether the edges an their terminal vertices are in the same orbit or not.

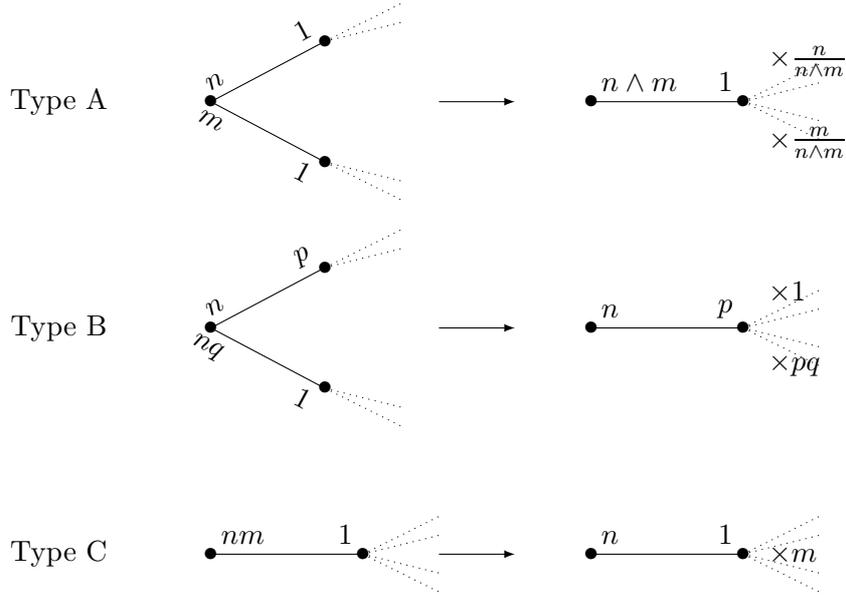
\begin{figure}
 \centering
 \begin{tikzpicture}
\begin{scope}
\draw (0,0) node {$\bullet$}  -- (1.5,0.8) node [near start, sloped, above left] {$n$} node {$\bullet$} node [near end, sloped, above right] {$1$} ;
\draw (0,0)  -- (1.5,-0.8)node [near start, sloped, below left] {$m$}  node [near end, sloped, below right] {$1$}node {$\bullet$};
\draw[dotted] (1.5,0.8) -- +(1, 0.5);
\draw[dotted] (1.5,0.8) -- +(1, 0.25);
\draw[dotted] (1.5,-0.8) -- +(1, -0.25);
\draw[dotted] (1.5,-0.8) -- +(1, -0.5);

\draw[-latex] (3,0)--(4,0);
\draw (5,0) node [above right] {$n \wedge m$} node {$\bullet$} -- (7,0) node [above left] {$1$} node {$\bullet$} ;
\draw[dotted] (7,0) node [above right = 0.2cm] {$\times \frac{n}{n\wedge m}$} -- +(1, 0.5);
\draw[dotted] (7,0) -- +(1, 0.25);
\draw[dotted] (7,0) -- +(1, -0.25);
\draw[dotted] (7,0) node [below right  = 0.2cm] {$\times \frac{m}{n\wedge m}$}-- +(1, -0.5);

\draw (-2,0) node {Type A};
\end{scope}
\begin{scope}[yshift = -3cm]
\draw (0,0) node {$\bullet$}  -- (1.5,0.8) node [near start, sloped, above left] {$n$} node {$\bullet$} node [near end, sloped, above right] {$p$} ;
\draw (0,0)   -- (1.5,-0.8) node [near start, sloped, below left] {$nq$} node [near end, sloped, below right] {$1$}node {$\bullet$};
\draw[dotted] (1.5,0.8) -- +(1, 0.5);
\draw[dotted] (1.5,0.8) -- +(1, 0.25);
\draw[dotted] (1.5,-0.8) -- +(1, -0.25);
\draw[dotted] (1.5,-0.8) -- +(1, -0.5);

\draw[-latex] (3,0)--(4,0);
\draw (5,0) node [above right] {$n$} node {$\bullet$} -- (7,0) node [above left] {$p$} node {$\bullet$} ;
\draw[dotted] (7,0) node [above right = 0.2cm] {$\times 1$} -- +(1, 0.5);
\draw[dotted] (7,0) -- +(1, 0.25);
\draw[dotted] (7,0) -- +(1, -0.25);
\draw[dotted] (7,0) node [below right  = 0.2cm] {$\times pq$}-- +(1, -0.5);

\draw (-2,0) node {Type B};
\end{scope}

\begin{scope}[yshift=-6cm]
\draw (0,0) node [above right] {$nm$} node {$\bullet$} -- (2,0) node [above left] {$1$}node {$\bullet$};
\draw[dotted] (2,0) -- +(1, 0.5);
\draw[dotted] (2,0) -- +(1, 0.25);
\draw[dotted] (2,0) -- +(1, -0.25);
\draw[dotted] (2,0) -- +(1, -0.5);

\draw[-latex] (3,0)--(4,0);
\draw (5,0) node [above right] {$n$} node {$\bullet$} -- (7,0) node [above left] {$1$}node {$\bullet$};
\draw[dotted] (7,0) node [right = 0.2cm] {$\times m$} -- +(1, 0.5);
\draw[dotted] (7,0) -- +(1, 0.25);
\draw[dotted] (7,0) -- +(1, -0.25);
\draw[dotted] (7,0) -- +(1, -0.5);

\draw (-2,0) node {Type C};
\end{scope}
\end{tikzpicture}
 \caption{The three types of folds which preserve the deformation space, as seen in labelled graphs.} \label{fig:typesPliagesDef}
\end{figure}

Denote by $A=G_v$, $B_i = G_{t(e_i)}$, $E_i = G_{e_i}$. 
The result of a fold is in the same deformation space as the original tree if and only if one of the following conditions\footnote{If we work with restricted deformations spaces $\D^{\mathcal A}$, we must ensure when folding that the new edge groups belong to $\mathcal{A}$. However in this article we will work backwards: starting with the folded tree, we will unfold it, and we do not need to worry since the new groups will be smaller.} is true (see figure \ref{fig:typesPliagesDef}):

\begin{itemize}
 \item the fold is of type I with $B_1=E_1$ and $B_2=E_2$, which we call \emph{type A} (not referring to subtypes of \cite{BestvinaFeighnBounding})
 \item the fold is of type I and $B_2 \subset B_1$ or $B_1 \subset B_2$ which we call \emph{type B}
 \item the fold is of type II and $E_1=B_1$ which we call \emph{type C}.
\end{itemize}

\begin{rema} \label{rem:etoile}
The pre-image of any edge by a fold which does not change the deformation space is star-shaped, i.e. it consists in a collection of edges which share a common endpoint. In fact 
\begin{itemize}
 \item $\pi ^{-1} \pi(\epsilon) = \epsilon$ if $\epsilon \notin G\cdot {e_1,e_2}$
 \item $\pi ^{-1} \pi(e_1) = G_{\pi(e_1)} \cdot \{e_1 \cup e_2\}$ and $G_{\pi(e_1)}$ fixes $o(e_1) = o(e_2)$ in case A, B and C.
\end{itemize}

\end{rema}

An \emph{unfolding of type C} is the inverse of a fold of type C, i.e. given a tree $T_0$, it consists in finding a tree $T_1$ such that $T_1 \to T_0$ is a fold of type C. We will need the following result about sequences of unfoldings. A \emph{sequence of unfoldings} of type C is a sequence $\dots \to T_n \to \dots \to T_1 \to T_0$ such that every $T_{i+1} \to T_i$ is a fold of type C.

\begin{lem} \label{lem:pasdexplosion}
 Suppose that $G$ is not a solvable Baumslag-Solitar group. Let $T$ be a GBS tree for $G$. Any sequence of unfoldings of type C (see Figure \ref{fig:typesPliages}) of $T$ is finite. 
\end{lem}

\begin{proof}
 We will show that it is impossible to construct an infinite sequence of unfoldings of type C. We consider the quotient $\Gamma = T/G$. Unfoldings of type C do not change the edges of $\Gamma$: only labels vary. Therefore in the rest of the proof edges will keep the same name after unfolding. Consider the product  $\displaystyle  \prod _{e \in T/G} |\lambda(e)|$ of all labels in $\Gamma$ (see Figure \ref{fig:produitDecroissant}). It is a positive integer.

 An unfolding of type C can be performed on an edge $e \in \Gamma$ if and only if $\lambda(\bar e)= \pm 1$ and there exists an integer $q$ such that $|q| > 1$ and $q$ divides all other labels at $t(e)$.

During the unfolding, the labels at only two vertices may change : $\lambda (e'), e' \in \E_{t(e)} \setminus \{e\}$ are divided by $q$ whereas $\lambda(e)$ is multiplied by $q$. Let $k$ be the valence of the vertex $t(e)$, where $k \geq 2$ by minimality. The product $\displaystyle  \prod _{e \in T/G} |\lambda(e)|$ is multiplied by $q^{2-k}$. If the valence of $t(e)$ in $\Gamma$ is at least 3 (i.e. $k\geq 3$) then the product of labels decreases. An unfolding at a valence 2 vertex does not change the product. To conclude we need to show that one cannot produce an infinite sequence of unfoldings of type C on edges whose terminal vertex has valence 2.

 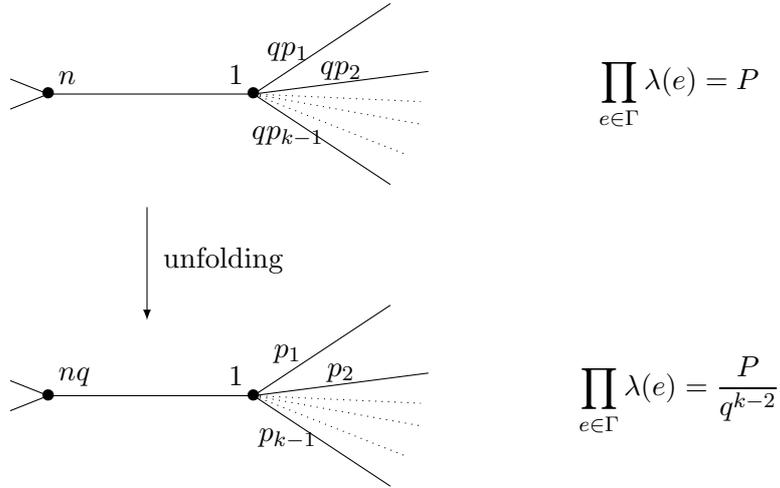
\begin{figure}
 \centering
 \begin{tikzpicture}[]
%\draw [dotted] (0,-6) rectangle (11.5,6);

%le haut
\begin{scope}[]

	\coordinate (w) at (0.7,0);
	\coordinate (v) at (3.4,0);

	%arete milieu
	\draw (w) node {$\bullet$} node [above right] {$n$} -- (v) node {$\bullet$} node [above left] {$1$};

	%arêtes droite
	\draw (v)  -- +(1.8,1.2)  node [near start, above] {$qp_1$};
	\draw (v)  -- +(2.3,0.3)  node [midway, above=-0.1cm] {$qp_2$};
	\draw (v)  -- +(1.8,-1.2)  node [near start, below] {$qp_{k-1}$};
	\draw[dotted] (v)  -- +(2.2,-0.1) ;
	\draw[dotted] (v)  -- +(2.2,-0.4); 
	\draw[dotted] (v)  -- +(2,-0.8); 
	%arêtes gauche
	\draw (w) -- +(-0.5,0.2);
	\draw (w) -- +(-0.5,-0.2);

	\node (pi) at (9,0) {$\displaystyle \prod_{e\in \Gamma} \lambda(e) = P$};
\end{scope}
%la flèche
	\draw[-latex] (2,-1.5) -- (2,-3);
	\node (leg) at (3, -2.2)  {unfolding};
% dessin du bas 
\begin{scope}[yshift=-4cm]

	\coordinate (w) at (0.7,0);
	\coordinate (v) at (3.4,0);

	%arete milieu
	\draw (w) node {$\bullet$} node [above right] {$nq$} -- (v) node {$\bullet$} node [above left] {$1$};

	%arêtes droite
	\draw (v)  -- +(1.8,1.2)  node [near start, above] {$p_1$};
	\draw (v)  -- +(2.3,0.3)  node [midway, above=-0.1cm] {$p_2$};
	\draw (v)  -- +(1.8,-1.2)  node [near start, below] {$p_{k-1}$};
	\draw[dotted] (v)  -- +(2.2,-0.1) ;
	\draw[dotted] (v)  -- +(2.2,-0.4); 
	\draw[dotted] (v)  -- +(2,-0.8); 
	%arêtes gauche
	\draw (w) -- +(-0.5,0.2);
	\draw (w) -- +(-0.5,-0.2);

	\node (pi) at (9,0) {$\displaystyle \prod_{e\in \Gamma} \lambda(e) = \frac{P}{q^{k-2}}$};

\end{scope}

\end{tikzpicture}
 \caption{The product decreases when doing unfoldings of type C.} \label{fig:produitDecroissant}
\end{figure}

In the rest of this proof we will denote this sort of unfolding by ``unfolding from a valence 2 vertex''. We will say that the unfolding of an edge $e$ is an unfolding ``from $t(e)$''. Since no new edge is created in $\Gamma$, edges will keep the same name after unfolding.

The \emph{topological edges} of the graph $\Gamma$ are connected components of the graph without its vertices with valence greater or equal to 3.
Since the number of topological edges of the graph is finite any infinite sequence of unfoldings from valence 2 vertices would have infinitely many unfoldings in at least one topological edge.

Note that an unfolding from a valence 2 vertex does not change any label outside of its topological edge. Consequently valence 2 unfoldings in different topological edges commute.

In most cases a topological edge is a segment $c$. It may be a circle, only if $\Gamma$ is a circle. First let us study the case of a segment (Figure \ref{fig:areteTopo}). It is composed of a certain number of edges with labels, and is bounded by vertices of valence 1 or at least 3. We will show that there is a bound (depending on $\Gamma$) on the length of any sequence of unfoldings from inner vertices in the same topological edge. As unfoldings in different topological edges commute, this will show that no infinite sequence of unfoldings from valence 2 vertices exist.

 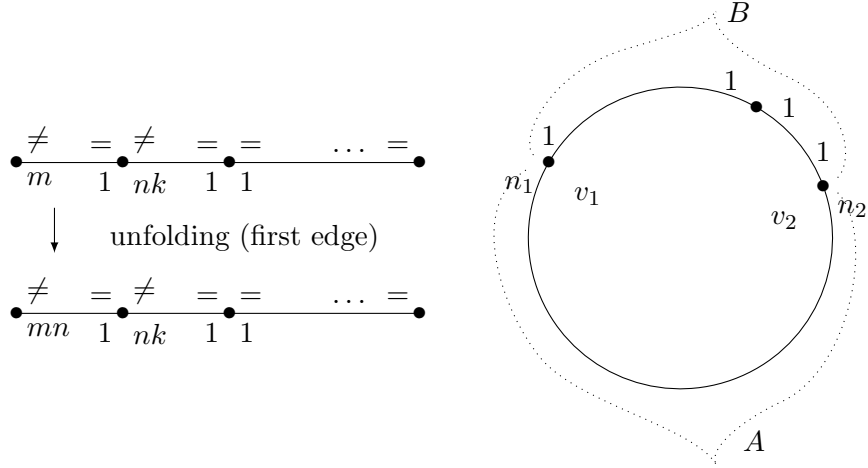
\begin{figure}
 \centering
 \begin{tikzpicture}[]
%\draw [dotted] (0,10) rectangle (12,-10);
\begin{scope}[]
	\draw (0,0) -- (5.3,0) ;

\coordinate (1) at (1.4,0);
\draw (1)  node {$\bullet$} node [above left] {$=$} node [above right] {$\neq$} ;
\draw (1)  node [below left] {$1$} node [below right] {$nk$} ;

\coordinate (2) at (2.8,0);
\draw (2)  node {$\bullet$} node [above left] {$=$} node [above right] {$=$} ;
\draw (2)  node [below left] {$1$} node [below right] {$1$} ;

\coordinate (3) at (4.8,0);
\draw (3)  node [above left] {$\dots$};

\coordinate (4) at (5.3,0);
\draw (4)  node {$\bullet$} node [above left] {$=$} ;

\coordinate (0) at (0,0);
\draw (0)  node {$\bullet$} node [above right] {$\neq$} ;
\draw (0)  node [below right] {$m$} ;

%flèche
\draw [-latex] (0.5,-0.6) -- (0.5,-1.2);
\node (dep) at (3, -1) {unfolding (first edge)};

	\draw (0,-2) -- (5.3,-2) ;

\coordinate (1) at (1.4,-2);
\draw (1)  node {$\bullet$} node [above left] {$=$} node [above right] {$\neq$} ;
\draw (1)  node [below left] {$1$} node [below right] {$nk$} ;

\coordinate (2) at (2.8,-2);
\draw (2)  node {$\bullet$} node [above left] {$=$} node [above right] {$=$} ;
\draw (2)  node [below left] {$1$} node [below right] {$1$} ;

\coordinate (3) at (4.8,-2);
\draw (3)  node [above left] {$\dots$};

\coordinate (4) at (5.3,-2);
\draw (4)  node {$\bullet$} node [above left] {$=$} ;

\coordinate (0) at (0,-2);
\draw (0)  node {$\bullet$} node [above right] {$\neq$} ;
\draw (0)  node [below right] {$mn$} ;
\end{scope}

\begin{scope}[xshift=9cm,yshift=-1cm]
% 	 \coordinate (c) at (0,-0);
% 	\draw (c) circle [radius=2cm];
% 	\node (d1) at  (30/17,16/17) {$\bullet$};
% 	\node (d2) at  (-30/17,16/17) {$\bullet$};
% 	\node (d3) at  (6/5,8/5) {$\bullet$};

	\coordinate (d1) at (-2,1);
	\draw (d1)  node {$\bullet$} arc (150:60:2) node (d2) {$\bullet$}  arc (60:20:2) node (d3){$\bullet$} arc (20:-210:2) ;
	\draw (d1) node [below right=0.2cm] {$v_1$} node [below left =0.05cm] {$n_1$} node [above=0.1cm] {$1$};
	\draw (d3) node [below left=0.2cm] {$v_2$} node [below right =0.05cm] {$n_2$} node [above = 0.2cm] {$1$};
	\draw (d2) node [right=0.2cm] {$1$} node [above left=0.1cm] {$1$};
	%courbes A et B
	\draw [dotted] (-2.2,1.1) .. controls +(-0.2,0.2) and (-2,2.2) .. (-1,2.5) .. controls (0,2.8) and (0,2.8) .. (0.2,3);
	\draw [dotted] (1.8,0.8) .. controls +(0.2,0.2) and (2,1.8) .. (1,2.2) .. controls (0.5,2.45) and (0.1,2.8) .. (0.2,3);
	\node (B) at (0.5,3) {$B$};

	\draw [dotted] (-2.3,0.9) .. controls +(-0.8,-0.5) and (-2.6,-1) .. (-2.1,-1.6)  .. controls (-1.6,-2.2) and (-0.1,-2.5) .. (0.2,-3);
	\draw [dotted] (1.8,0.6) .. controls +(0.5,-0.5) and (2,-1.8) .. (1,-2.2) .. controls (0.5,-2.45) and (0.1,-2.8) .. (0.2,-3);
	\node (B) at (0.7,-2.7) {$A$};

\end{scope}
\end{tikzpicture}
 \caption{Topological edges} \label{fig:areteTopo}
\end{figure}
 
 Choose an orientation for $c$ and view it as a concatenation of oriented edges $e_1 \dots e_l$. We define the following complexity for a segment $c$ of length $l$:
  \[
  K(c):= \prod _{i = 1} ^n |\lambda(e_i)|^i \times \prod _{i=1}^n |\lambda(\bar e_i)|^{n+1-i}
 \]
 which gives more weight to initial labels of edges at the end of the chain and to terminal labels of edges at the beginning of the chain. Note that it is well-defined since it does not depend on the choice of orientation for $c$.
 
 The idea is that unfoldings move factors of labels in the direction where they will weight less. Let us show that the complexity decreases during any unfolding.

 Let $i \in \{1, \dots, l-1\}$. Suppose we do an unfolding of $e_i$ (a valence 2 unfolding, so we exclude the case $i=l$). We obtain a new segment $c'$ with different labels. Two labels change:  $\lambda(e_{i+1})$ is divided by a factor $q$ and $\lambda(e_i)$ is multiplied by this same factor $q$, where $|q| \neq 1$. Therefore $K(c') = K(c) \times |q|^i /{|q|^{i+1}} = K(c) / |q|$ so $K(c') < K(c)$.

 Since the formula for the complexity does not depend on the orientation of the segment we get the same result for an unfolding of a $\bar e_i$.

 As $K$ is a positive integer it cannot decrease indefinitely. This proves that the sequence of unfoldings must stop when the topological edge is a segment.

  \bigskip
  
 When the graph is a circle, it has a single topological edge. If there exists a vertex group which is strictly bigger than all edge groups at this vertex, then no unfolding can occur from this vertex and labels at this vertex can only increase. The graph minus this vertex is a segment to which we can apply the argument above: no unfolding sequence can be infinite.
  
 If not then every vertex has at least one label which is $\pm 1$. Choose an orientation of the circle and orient the edges accordingly. If all initial labels (for this orientation) are $\pm 1$ then the group is a solvable Baumslag-Solitar group. The same deduction can be made with the reverse orientation. If this happens for neither orientation then we can find two vertices with a label different from $\pm 1$ pointing in different directions. Call $v_1, v_2$ these two vertices, which split the circle into two segments. Call $B$ the segment with the $\pm 1$ labels at its endpoints and $A$ the segment with greater labels. Up to taking a subsegment of $B$ we may suppose that the labels borne by the edges of $B$ are all $\pm 1$.
 
 Unfoldings in $A$ cannot be done indefinitely, because no unfolding from $v_1$ or $v_2$ occur and we are again in the case of a segment. This shows that any long enough sequence of unfoldings must involve an edge in $B$.
  
 The only unfoldings which may involve edges in $B$ are unfoldings from $v_1$ or $v_2$. If we perform such an unfolding we increase one of the labels inside $B$.
 
 Then we can define a new partition where $B$ strictly decreases, and iterate until there is a vertex with both labels different from $\pm 1$. 
\end{proof}

\section{Special factors, and algorithm for simplicity of group elements} \label{sec:facteursSpeciaux}

\subsection{Special factors for a deformation space}

In this section, we introduce special factors with respect to some deformation space $\D^{\A}$. They are an analogue of free factors for free groups.

\begin{rema}
 We will write $\bar \D^{\A}$ to denote the set of all trees obtained by collapsing $G$-invariant subforests in trees of $\D^{\A}$, including the trivial tree. We include trees of $\D^{\A}$ which correspond to collapsing empty forests.
\end{rema}

\begin{defi}
 A \emph{special factor $H$ with respect to $\D^\A$} is a subgroup of $G$ which is the stabilizer of a point in a tree $T$ in $\bar \D^\A$ and which is not elliptic in $\D^\A$.
  
 When the deformation space is obvious, we will write simply \emph{special factor}. 
 
 We call $H$ a \emph{proper special factor} when $H \neq G$. Elliptic groups, i.e. vertex stabilizers of trees in $\D^\A$, are not considered to be actual special factors.
\end{defi}

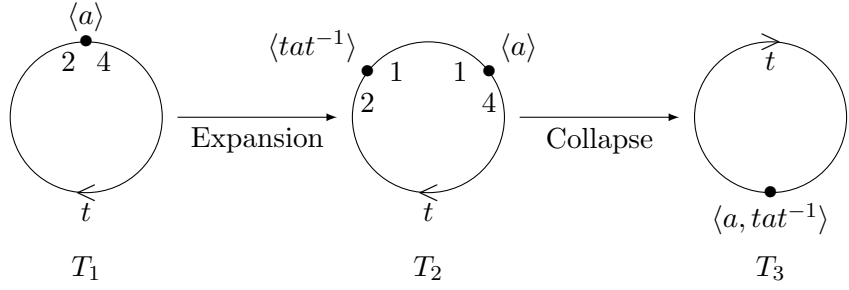
\begin{figure}
 \centering
 \begin{tikzpicture}[]
%cercle de gauche
\draw (-0.5,1) arc (90:450:1)node [midway] {$<$}  node [midway, below] {$t$};
\draw (-0.5,1) node {$\bullet$} node [above] {$\langle a \rangle$} node [below left] {$2$} node [below right] {$4$};

\draw (-0.5, -2) node {$T_1$};

% flèche
\draw[-latex] (0.7,0) -- (2.8,0) node [midway,below] {Expansion};

%cercle du milieu
\draw (4,1) arc (90:450:1)node [midway] {$<$} node [midway, below] {$t$};
\draw (4.8,0.6) node {$\bullet$} node [above right] {$\langle a \rangle$} node 
[ left = 0.15cm] {$1$} node [below = 0.15cm] {$4$};
\draw (3.2,0.6) node {$\bullet$} node [above left] {$\langle t a t^{-1} \rangle$} node 
[ right = 0.15cm] {$1$} node [below = 0.15cm] {$2$};

\draw (4, -2) node {$T_2$};

% flèche de droite
\draw[-latex] (5.2,0) -- (7.3,0) node [midway,below] {Collapse};

% cercle de droite
\draw (8.5,-1) arc (-90:270:1)node [midway] {$>$}  node [midway, below] {$t$};% node [left=0.7cm, rotate=-50] {$ \langle t a t^{-1} \rangle \rightarrow$} node [right=0.7cm, rotate=50] {$\leftarrow \langle a \rangle$};
\draw (8.5,-1) node {$\bullet$} node [below] {$\langle a, tat^{-1} \rangle $};

\draw (8.5, -2) node {$T_3$};

\end{tikzpicture}
 \caption{Expansion and collapse on the standard tree $T_1$ (seen in the quotient) and new special factor with respect to $\D$ obtained after the collapse, visible in $T_3 \in \bar \D$. The tree $T_2$ is in $\D$ but not in $\D^{\Amin}$.} \label{fig:facteursSpeciaux}
\end{figure}

\begin{ex}
 This notion depends on the allowed edge groups in $\D^{\A}$. 
 The space $\D^{\Amin}$ of reduced trees has fewer allowed edge groups than $\D$. Consider the standard tree $T_1$ for $BS(2,4):=\langle a,t | ta^2t^{-1}=a^4 \rangle$ (see Figure \ref{fig:facteursSpeciaux}) and perform an expansion (yielding $T_2$) and a collapse (yielding $T_3$) as described by the figure. We obtain a special factor with respect to $\D$, which is the subgroup $\langle a, tat^{-1}\rangle$. This subgroup cannot be obtained by collapsing a tree in $\D^{\Amin}$, thus it is not a special factor with respect to $\D^{\Amin}$. 
\end{ex}

For a subgroup $H <G$, denote its conjugacy class by $[H]$.

\begin{defi} \label{defi:systeme}
 A \emph{system of special factors with respect to $\D^\A$} is a finite collection of conjugacy classes of subgroups $\mathcal H := \{ [H_1], \dots, [H_k]\}$ of $G$ such that there exists $T_{\mathcal H} \in \bar \D^\A$ such that $\mathcal H$ is the set of conjugacy classes of vertex stabilizers in $T_{\mathcal H}$ which are not elliptic in $\D$.
 
 The system is \emph{proper} if it is not $\{[G]\}$.
\end{defi}
\begin{rema}
 Just like special factors, a system of special factors can be viewed in a graph of groups. Is is given by a collection of disjoint subgraphs $\Gamma_1, \dots, \Gamma_k$ of $\Gamma$ such that for every $i \in \{1, \dots, k\}$ the subgroup $H_i$ is isomorphic to $\pi_1(\Gamma_i)$.
\end{rema}

\begin{defi}
  Let $\mathcal H := \{[H_1], \dots, [H_k]\}$ be a system of special factors.
  We say that a collection $\G$ of elements of $G$ is \emph{$\mathcal H$-peripheral}, which we write $\G \preceq \mathcal H$, if for any $g \in \G$ there exists $1 \leq i \leq k$ such that $g$ is contained in a conjugate of $H_i$.
  
  The collection $\G$ is \emph{simple} if there exists a proper system of special factors $\mathcal H$ such that $\G \preceq \mathcal H$.
 
 A system of cyclic factors $\mathcal H'$ is \emph{$\mathcal H$-peripheral} ($\mathcal H' \preceq \mathcal H$) if for every conjugacy class $[H']\in \mathcal H'$ there exists $[H] \in \mathcal H$ such that $H'$ can be conjugated into a subgroup of $H$. 
\end{defi}
\begin{remas}\label{rem:ordre}
 \begin{enumerate}
  \item Equivalently $\mathcal H' \preceq \mathcal H$ if there exists a $G$-equivariant map $T_{\mathcal H'} \to T_{\mathcal H}$ where $T_{\mathcal H'}, T_{\mathcal H}$ are defined as in Definition \ref{defi:systeme}.
  \item The relation $\preceq$ defines an order on the set of systems of special factors. It is obviously reflexive and transitive. Suppose $\mathcal H' \preceq \mathcal H$ and $\mathcal H \preceq \mathcal H'$, then we get two maps $T_{\mathcal H'} \to T_{\mathcal H}$ and $T_{\mathcal H} \to T_{\mathcal H'}$. This implies that $T_{\mathcal H'}$ and $T_{\mathcal H}$ have the same elliptic subgroups. In particular they have the same non-cyclic vertex stabilizers, so $\mathcal H = \mathcal H'$. Thus $\preceq$ is antisymmetric. 
 \end{enumerate}
\end{remas}

\bigskip

\subsection{Whitehead graph and criterion for simplicity}

We fix a collection of cyclic allowed edge groups $\mathcal A$ and consider a restricted deformation space $\D^\A$.

\begin{defi}
 Let $T$ be a $G$-tree. A \emph{turn} in $T$ is an unordered pair of distinct edges with same origin. If $e,e'$ are two such edges, the corresponding turn is denoted by $\{e, e'\}$.
 
 When $e=e'$ we call the pair a \emph{degenerate turn}. 

 A geodesic $\gamma$ \emph{crosses} a turn $\{e,e'\}$ if $\gamma$ contains $e \cup e'$.
\end{defi}

Let $\mathcal G :=\{g_1, \dots, g_k\} \in G$ be a finite collection of loxodromic elements. Let $T \in \D^\A$ and $v$ a vertex in $T$. The set $\E _v$ is the set of edges of $T$ with origin $v$.

\begin{defi}
 The \emph{Whitehead graph} $Wh_T(\G,v)$ is the following graph. The vertex set is $\E_v$. Two vertices $e,e'$ are linked by a non-oriented edge in $Wh_T(\G,v)$ when there exists $g \in \G$ and some conjugate of $g$ whose axis crosses the turn $\{e, e'\}$.
\end{defi}

\begin{remas}
 \begin{enumerate}
  \item Equivalently, we link $e$ and $e'$ by an edge whenever there is $h \in G $ such that $\{he, he'\}$ is a turn in the axis of some $g \in \G$.
  \item The Whitehead graph is a simplicial graph. In particular it does not have any loop.
  \item The group $G_v$ has a natural action on $Wh_T(\G,v)$.
 \end{enumerate}
\end{remas}

\begin{defi}
 When $Wh_T(\G,v)$ is not connected, we call \emph{admissible connected component} any connected component in $Wh_T(\G,v)$ whose stabilizer is in $\A$. When all edge groups are allowed all connected components are automatically admissible.
 
 Let $p$ be a vertex in $Wh_T(\G,v)$. Let $W_0$ be the connected component which contains $p$. The vertex $p$ is an \emph{admissible cut point} if $W_0 \setminus \{p\}$ is disconnected and if there exists a connected component $A$ of $W_0 \setminus \{p\}$ satisfying $A \cap G_v \cdot p = \varnothing$.
 
 The Whitehead graph $Wh_T(\G,v)$ has an \emph{admissible cut} (for $\A$) when it has either an admissible connected component or an admissible cut point.
\end{defi}

\begin{remas} \label{rem:stabilisateur}
 \begin{enumerate}
  \item  In the admissible cut point definition, the stabilizer of $A$ is automatically an allowed edge group since it is a subgroup of $G_{e_p}$, where $e_p$ is the edge of $T$ corresponding to the vertex $p$ of the Whitehead graph.
  \item Since the Whitehead graph has no loop, $A$ contains a vertex.
 \end{enumerate}
\end{remas}

The following lemma uses that $T$ is locally finite in an essential way.

\begin{lem}[Dual tree to the Whitehead graph] \label{lem:arbreDual}
  Let $p$ a vertex in $Wh_T(\G,v)$ and $W_0$ the connected component containing $p$. If $p$ is a cut point in $W_0$ (i.e. $W_0 \setminus \{p\}$ is not connected) then $p$ is an admissible cut point of the Whitehead graph.
\end{lem}

\begin{proof}
 Let $p$ a vertex in $Wh_T(\G,v)$ whose complement is disconnected. A dual forest to the Whitehead graph can be defined as follows.

 First we define the following equivalence relation on the geometric realization of the Whitehead graph $W$: for all points $x,y \in W$, $x \sim y$ if for all $q \in G_v \cdot p$, $x$ and $y$ are in the same connected component of $W \setminus \{q\}$. Equivalence classes of this relation define a partition of $W \setminus G_v \cdot p$. An equivalence class may contain no vertex of $W$ and that is why we work with the geometric realization of the graph. This partition is coarser than the partition into connected components of $W \setminus G_v \cdot p$.

 Then we define the bipartite graph $B$ as follows. There is a vertex $u_q$ for every vertex $q \in G_v \cdot p$. There is also a vertex $v_P$ for every equivalence class $P$ of the equivalence relation defined above. We put an edge between $u_q$ and $v_P$ if $q\in \bar P$.

 The graph $B$ obtained is a forest because every vertex $u_q$ disconnects all its neighbours in $B$. It is connected if and only if the Whitehead graph is.
 
 \medskip
 
 Let $W_0$ be the connected component of $W$ containing $p$. Suppose $p$ is a cut point of $W_0$.
 The component $B_0$ of $B$ containing $u_p$ is a finite tree so it has a terminal vertex $w$. 
 Vertices $u_q$ cannot be terminal since $q$ is a cut vertex so $w = v_P$ for some equivalence class $P$, which is attached to a vertex $u_q=h\cdot u_p$ of $B_0$.
 
 The equivalence class $P$ is a connected component of $W_0 \setminus \{h \cdot p\}$ because $G_v \cdot u_p \cap \bar P = \{u_q\}$. 
 Taking $A=P$ in the definition, $h \cdot p$ (and thus $p$) is an admissible cut point.
\end{proof}

We can now state the main theorem of this section. Its proof is given in subsection \ref{subsec:proofOfTheo}.
 
\begin{theo} \label{theo:Whitehead}
 Let $\G \in G$ be a finite collection of loxodromic elements. If $\G$ is simple with respect to $\D^\A$ then for all $T \in \D^\A$ there exists a vertex $v \in T$ such that $Wh_T(\G,v)$ has an admissible cut for $\A$.
\end{theo}

\begin{rema}
 A solvable Baumslag-Solitar group $BS(1,n):=\langle a,t | tat^{-1}=a^n \rangle$ has no proper special factor (with respect to $\D=\D^{\Amin}$). 
\end{rema}

\begin{rema}\label{rem:subdivision}
 Let $T$ a $G$-tree and $T'$ a tree obtained from $T$ by subdividing an edge. Then there exists a vertex $v \in T$ such that $Wh_T(\G,v)$ has an admissible cut if and only if there exists such a vertex in $T'$.

 Indeed $T'$ inherits the Whitehead graphs of $T$ in addition with another Whitehead graph coming from the additional vertex $v'$. The latter is a graph containing exactly two vertices. If they are joined by an edge, there is no admissible cut. If not, $Wh_{T'}(\G,v')$ is disconnected, which means for any $g \in \G$, no translate of the axis of $g$ crosses the subdivided edge. In $T$ this edge must then appear as an isolated vertex in the Whitehead graph of one of its endpoints, so some Whitehead graph in $T$ has an admissible cut.
\end{rema}

\bigskip

\subsection{Unfolding lemma}

 Let $T \in \D^\A$. According to remark \ref{rem:subdivision} we may assume the following : up to performing a finite number of edge subdivisions at the beginning, $T$ has no edge with both ends in the same orbit, i.e. $T/G$ has no loop. This allows us to deal with fewer cases in the proof. The proof is similar to that of an analogous result concerning the case of free products in \cite[Proposition 5.1]{GuirardelHorbezAlgebraic}. We will need the following lemma which enables us to perform unfoldings on $T$ or expansions when we find a Whitehead graph with an admissible cut. We allow vertices of valence $2$ in the trees considered.

\begin{lem} \label{lem:grapheWhEtDepliage}
 Suppose $T/G$ has no loop.
 The following conditions are equivalent :
 \begin{enumerate}
  \item There exists a Whitehead graph $Wh_T(\G,v)$ with an admissible cut with respect to $\D^\A$.
  \item There exists a tree $S \in \D^\A$ and a non-injective $G$-equivariant application $f: S \to T$ sending edge to edge or edge to vertex such that for every $g \in \G$, $\|g\|_S=\|g\|_T$. 
 \end{enumerate}
 In the second condition $S/G$ has also no loop. Moreover the map $f$ can be chosen to be either a fold or a collapse.
\end{lem}

Let us start with a preliminary result about lifting the axis of an element $g \in \G$ when performing an unfolding.

\begin{lem} \label{lem:relevementLocal}
 Let $g$ be a loxodromic element in $G$.
 Let $T, S \in \D^\A$ and $f: S \to T$ a simplicial map such that for all edge $e \in T$, edges in the pre-image $f^{-1}(e):=\{ \tilde e \in E(S) / f(\tilde e)= e\}$ all share a common vertex. Suppose that every turn in the axis of $g$ lifts to $S$, that is to say: for every turn $\{e_1,e_2\}$ in $\axe_T(g)$ there exists a turn $\{\tilde e_1, \tilde e_2\}$ in $S$ such that $f(\tilde e_1) = e_1, f(\tilde e_2) = e_2$.
 
 Then $\axe_T(g)$ lifts isometrically in $S$. Equivalently $f$ is isometric on $\axe_S(g)$.
\end{lem}

\begin{rema}
 Remark \ref{rem:etoile} states that this lemma applies to 
 folds of type A, B and C.
\end{rema}

\begin{figure}
  \centering
  \newcommand{\etoilecf}[3]{%
	\begin{scope}[shift={#1}, rotate={#3},scale=1.3]
	\draw[#2] (0,0) -- +(-1,0) node [midway,sloped] {$<$};
	%\draw (#2,2) node {$+$};
	\draw[#2] (0,0) -- +(1,0) node [midway,sloped] {$>$};
	\draw[#2] (0,0) -- +(-0.6,0.8) node [midway,sloped] {$<$};
	\draw[#2] (0,0) -- +(-0.6,-0.8) node [midway,sloped] {$<$};
	\draw[#2] (0,0) -- +(0.6,0.8) node [midway,sloped] {$>$};
	\draw[#2] (0,0) -- +(0.6,-0.8) node [midway,sloped] {$>$};
	\end{scope}
}
\newcommand{\etoilecp}[3]{%
	\begin{scope}[shift={#1}, rotate={#3},scale=1.3]
	\draw[#2] (0,0) -- +(-1,0) node [midway,sloped] {$>$};
	%\draw (#2,2) node {$+$};
	\draw[#2] (0,0) -- +(1,0) node [midway,sloped] {$<$};
	\draw[#2] (0,0) -- +(-0.6,0.8) node [midway,sloped] {$>$};
	\draw[#2] (0,0) -- +(-0.6,-0.8) node [midway,sloped] {$>$};
	\draw[#2] (0,0) -- +(0.6,0.8) node [midway,sloped] {$<$};
	\draw[#2] (0,0) -- +(0.6,-0.8) node [midway,sloped] {$<$};
	\end{scope}
}
\begin{tikzpicture}[]

\begin{scope}[xshift=0,yshift=0]
%E1 centrifuge, E2 centrifuge
\etoilecf{(0,0)}{blue}{0}
\etoilecf{(1.3,0)}{red}{10}
\draw [blue] (-0.3,1.3)  node {$A_1$};
\draw [red] (1.2,1.3)  node {$A_2$};

\draw (0.5,-1.3) node {$A_1$ and $A_2$ centrifugal};
\end{scope}

\begin{scope}[xshift=6cm,yshift=0]
%E1 centrifuge, E2 centripète
\etoilecf{(0,0)}{blue}{0}
\etoilecp{(2.57,0.21)}{red}{10}
\draw [blue] (-0.3,1.3)  node {$A_1$};
\draw [red] (2.3,1.3)  node {$A_2$};

\draw (1.2,-1.3) node {$A_1$ centrifugal and $A_2$ centripetal};
\end{scope}

\begin{scope}[xshift=0,yshift=-4cm]
%E1 centripète, E2 centrifuge
\etoilecp{(0.7,0)}{blue}{-15}
\etoilecf{(0.7,0)}{red}{15}
\draw [blue] (-0.3,1.3)  node {$A_1$};
\draw [red] (1.8,1.3)  node {$A_2$};

\draw (1,-1.5) node {$A_1$ centripetal and $A_2$ centrifugal};
\end{scope}

\begin{scope}[xshift=6cm,yshift=-4cm]
%E1 et E2 centripètes
\etoilecp{(0.5,0)}{blue}{0}
\etoilecp{(1.8,0.25)}{red}{10}
\draw [blue] (-0.3,1.3)  node {$A_1$};
\draw [red] (2.4,1.3)  node {$A_2$};

\draw (1.5,-1.5) node {$A_1$ and $A_2$ centripetal};
\end{scope}

\end{tikzpicture}
  \caption{Relative dispositions of the stars depending on their orientations} \label{fig:etoiles}
 \end{figure}
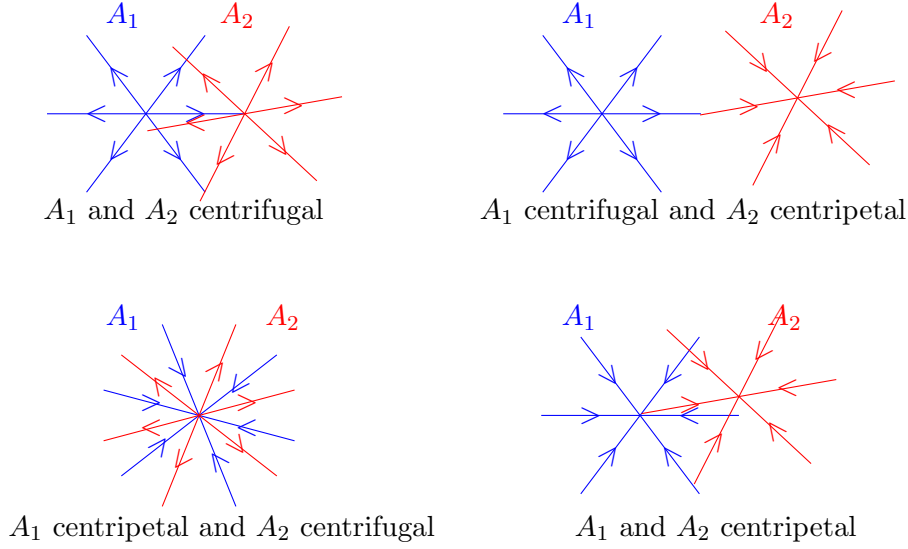

\begin{proof}
 First of all, given an orientation of an edge $e \in T$, the edges in $f^{-1}(e)$ get a compatible orientation. We will call a set of edges with a common vertex a \emph{star}. If $f^{-1}(e)$ is a star then this orientation is either centripetal or centrifugal.
 
 \medskip
 
 We claim that for every edge $e \in \axe_T(g)$ the intersection $f^{-1} (e) \cap \axe_S(g)$ consists in a unique edge $\tilde e$. Moreover, if $e, e'$ are adjacent in $T$ then $\tilde e, \tilde e'$ are adjacent in $S$. This yield a continuous application $\axe_T(g) \to \axe_S(g)$ which is an inverse for $f$ on the axis. This proves the lemma.
 
 \medskip 
 
 Now let us prove the claim. Let $e_1, e_2, e_3$ be three consecutive edges in $\axe_T(g)$ with $t(e_1)=o(e_2)$ and $t(e_2)=o(e_3)$.  We will show that $f^{-1}(e_2) \cap \axe_S(g)$ consists in exactly one edge.
 
 Let $A_1, A_2, A_3$ be the respective pre-images of $e_1, e_2, e_3$. They are stars. We endow them with an orientation, either centrifugal or centripetal, compatible with the orientation of their image. Because of orientations in $T$ and since the turns lift, the star $A_1$ is attached to an end of $A_2$ if the latter is centripetal, and to the centre if it is centrifugal (see figure \ref{fig:etoiles} for a picture of the different cases). On the contrary, $A_3$ is attached to the centre of $A_2$ if $A_2$ is centripetal and to an end if it is centrifugal. In both cases, distance between $A_1$ and $A_3$ is $1$. There is a unique edge in $A_2$ which is adjacent to both $A_1$ and $A_3$ and we call it $\tilde e_2$.
 
 Since $\axe_T(g) \subset f(\axe_S(g))$, $\axe_S(g)$ intersects both $A_1$ and $A_3$. Since it is a geodesic, its intersection with $A_2$ must be the single edge $\tilde e_2$. This proves the first part of the claim. 
 
 The second part follows: the lift $\tilde e_1$ is adjacent to $A_2$ and the intersection of two stars is a single point, so it is adjacent to $\tilde e_2$.
\end{proof}

\begin{proof} [Proof of lemma \ref{lem:grapheWhEtDepliage}]
 Suppose that the first condition is true: there is a vertex $v \in T$ such that $Wh_T(\G,v)$ has an admissible cut. We distinguish several cases based on the shape of $Wh_T(\G,v)$ and give an application $f:S \to T$ for each case  (see figure \ref{fig:casWhitehead}). The map $f$ will be a collapse in Case 1, a type A or type C fold in Case 2, and a type B fold in Case 3. Types of folds were defined in subsection \ref{subsec:folds}. Note that although type A and B folds may look similar, they lead to very different Whitehead graphs which need to be dealt with separately.
 
  \begin{figure}
  \centering
  %\usetikzlibrary{shapes}
\begin{tikzpicture}[scale=0.8]
	%Cas disconnexe (1)
	\begin{scope}[]
		\node (cas1) at (0.5,3) {Case 1: disconnected graph};
		
		\node[draw,ellipse] (A) at (0,0) {$A$};							        
		\node[draw,ellipse] (gA) at (0.5,1.5) {$gA$};
		\node[draw,ellipse] (B) at (2,0) {$B$};
		
		\draw (-2,1) node {$Wh_T(g,v)$};

		\draw[-latex] (1,-0.7) -- (1,-2);

		\draw (1,-4) node {$\bullet$} -- (0, -4.2) node {$\bullet$} ;
		\draw (1,-4) -- (0.1, -3.5) node {$\bullet$} ;
		\draw (1,-4) -- (2.2, -3.6) ;
		\draw (1,-4) -- (2.3, -3.9) node [right]{$B$} ;
		\draw (1,-4) -- (2.2, -4.2) ;
		%A
		\draw (0,-4.2) -- (-1.1,-4.1);
		\draw (0,-4.2) -- (-1.2,-4.3)  node [left]{$A$};
		\draw (0,-4.2) -- (-1.15,-4.5);
	
		\draw (-1,-3 ) node [left]{$gA$} -- (0.1, -3.5) ;
		\draw (-1.1,-3.2) -- (0.1, -3.5);
		\draw (-1.15,-3.35) -- (0.1, -3.5);
		
		\draw (-2.5, -3.5) node {$S$};
	\end{scope}

Cas 2 un seul pt coupure
	\begin{scope}[yshift=-9.5cm, xshift=2cm]
		\node[ centered] (cas1) at (3,3) {Case 2 : $G_v$-stable cutpoint};
		
		\node[draw,ellipse] (A) at (1,0) {$A$};							    
		\node[draw,ellipse] (gA) at (0.5,1.5) {$gA$};
		\node[draw,ellipse] (B) at (-1.1,0.3) {$B$};
		\coordinate (centre) at (0,0.5)  ;
		\draw (centre) node {$\bullet$} node [below] {$e$} --(B);
		\draw (centre) -- (A);	
		\draw (centre) -- (gA);
		\draw[-latex] (0,-0.7) -- (0,-2);
		
		\draw (-3,1.5) node {$Wh_T(g,v)$};

		\coordinate (racine) at (-2,-4);
		\draw (racine) node {$\bullet$} -- +(2,0.5) coordinate (point) node {$\bullet$} node [midway, sloped] {$>$} node [midway, above] {$e_A'$};
		\draw (racine)  -- +(1.8,-0.5) coordinate (autrepoint) node {$\bullet$} node [midway, sloped] {$>$} node [midway, below] {$e_B$};

		\draw (point)-- +(1,1) node [right]{$gA$};
		\draw (point)-- +(1.2,0.7);
		\draw (point)-- +(1.3,0.2) node [below right]{$A$};
		\draw (point)-- +(1.2,-0.2);

		\draw (autrepoint)-- +(1.2,0.4);
		\draw (autrepoint)-- +(1.3,0.0)  node [ right]{$B$};
		\draw (autrepoint) -- +(1.2,-0.4);
	
		\draw [dotted] (0.3,0.7) .. controls (-0.7,1) and  (-0.7,2.4) .. (0.5,2.5);
		\draw [dotted] (0.3,0.7) .. controls (0.7,0.45) and  (0.3,0.2) .. (0.3,0.2);
		\draw [dotted] (0.8,-0.8) .. controls (0.2,-0.6) and  (0.2,-0.1) .. (0.3,0.2);
		\draw [dotted] (0.8,-0.8) .. controls (2.8,-0.8) and  (1.7,2.5) .. (0.5,2.5) node [near end, above right] {$A'$};
        \draw (-3, -4.5) node {$S$};
%cas 2 bis
	\node (ou) at (3,0.8) {or};

		\node[draw,ellipse] (A) at (7,0) {$A$};							    
		\node[draw,ellipse] (gA) at (5.8,1.5) {$gA$};
		\node[draw,ellipse] (B) at (5.1,0) {$hA$};
		\coordinate (centre) at (6,0.5)  ;
		\draw (centre) node {$\bullet$} --(B);
		\draw (centre) -- (A);	
		\draw (centre) -- (gA);
		%\draw[-latex] (2,-0.7) -- (0,-2);
		
		\draw (8,1) node {$Wh_T(g,v)$};

		\draw[-latex] (6,-0.7) -- (6,-2);

		\coordinate (racine2) at (4.4,-3.6);
		\draw (racine2) node {$\bullet$} -- +(1.7,0.6) coordinate (point2) node {$\bullet$} node [midway, sloped] {$>$} node [midway, above] {$e'$};
		\draw (racine2)  -- +(1.9,-0.2) coordinate (autrepoint2) node {$\bullet$};
		\draw (racine2) -- +(1.6, -1) coordinate (derpoint) node {$\bullet$};

		%\draw (point2)-- +(1,1);
		\draw (point2)-- +(1.2,0.7) node [below right]{$A$};
		\draw (point2)-- +(1.3,0.2);
		%\draw (point2)-- +(1.2,-0.2);

		\draw (autrepoint2)-- +(1.2,0.3) node [below right]{$gA$};
		%\draw (autrepoint2)-- +(1.3,0.0);
		\draw (autrepoint2) -- +(1.2,-0.3);

		\draw (derpoint) -- +(1.3, 0.2)  node [below right]{$hA$};
		\draw (derpoint) -- +(1.2,-0.4);
		\draw (9, -4) node {$S$};

	\end{scope}

\begin{scope}[xshift=6.5cm]
		\node (cas1) at (1.5,3) {Case 3 : other kind of cutpoint};
		
		\node[draw,ellipse] (A) at (-0.8,-0.1) {$A$};
		\node[draw,ellipse] (g2A) at (0,1.5) {$g^2A$};
		\node[draw,ellipse] (B) at (1,0.2) {$B$};
		\node[draw,ellipse] (gA) at (3,-0.2) {$gA$};
		\node[draw,ellipse] (g3A) at (2.6,1.4) {$g^3A$};
		
		\draw (5,1) node {$Wh_T(g,v)$};

		\coordinate (centreg) at (0.1,0.35) ;
		\coordinate (centred) at (2,0.4);
		\draw (centreg) node [above left=0.05cm]{$e$} node {$\bullet$} -- (A) ;
		\draw (centreg) -- (g2A);
		\draw (centreg) -- (B);

		\draw (centred)  node [right=0.2cm]{$ge$} node {$\bullet$} -- (gA);
		\draw (centred) -- (g3A);
		\draw (centred) -- (B);

		\draw[-latex] (1,-0.7) -- (1,-2);

		\coordinate (v0) at (1.3,-3.7);
		\coordinate (w) at (-0.5,-3);
		\coordinate (gw) at (-0.5,-4);
		\coordinate (v) at (1,-2.8);
		\coordinate (gv) at (1,-4.6);

		\draw (v0) node {$\bullet$} -- (w) node {$\bullet$}node [midway, sloped] {$>$} node [midway, above right] {$e_1$};
		\draw (w) -- (v) node {$\bullet$} node [midway, sloped] {$>$} node [midway, above] {$e_2$};
		\draw (v0) -- (gw) node {$\bullet$};
		\draw (gw)-- (gv) node {$\bullet$};

		\draw (v0) -- +(1.2,0.3);
		\draw (v0) -- +(1.3,0.1) node [right] {$B$};
		\draw (v0) -- +(1.2,-0.2);

		\draw (v) -- +(1.2,0.8);
		\draw (v) -- +(1.3,0.6)node [ above right] {$g^2A$};
		\draw (v) -- +(1.4,0.2);
		\draw (v) -- +(1.3,-0.1) node [ right] {$A$};

		\draw (gv) -- +(1.4,0.1);
		\draw (gv) -- +(1.5,-0.1) node [right] {$g^3A$};
		\draw (gv) -- +(1.4,-0.5);
		\draw (gv) -- +(1.2,-0.7) node [right] {$gA$};
		
		\draw (4, -3.5) node {$S$};
\end{scope}

\end{tikzpicture}
  \caption{The three cases in the first part of the proof of the lemma. Above, the shape of the Whitehead graph; below, the shape of the tree around the corresponding vertex after transformation.}
 \label{fig:casWhitehead}
 \end{figure}
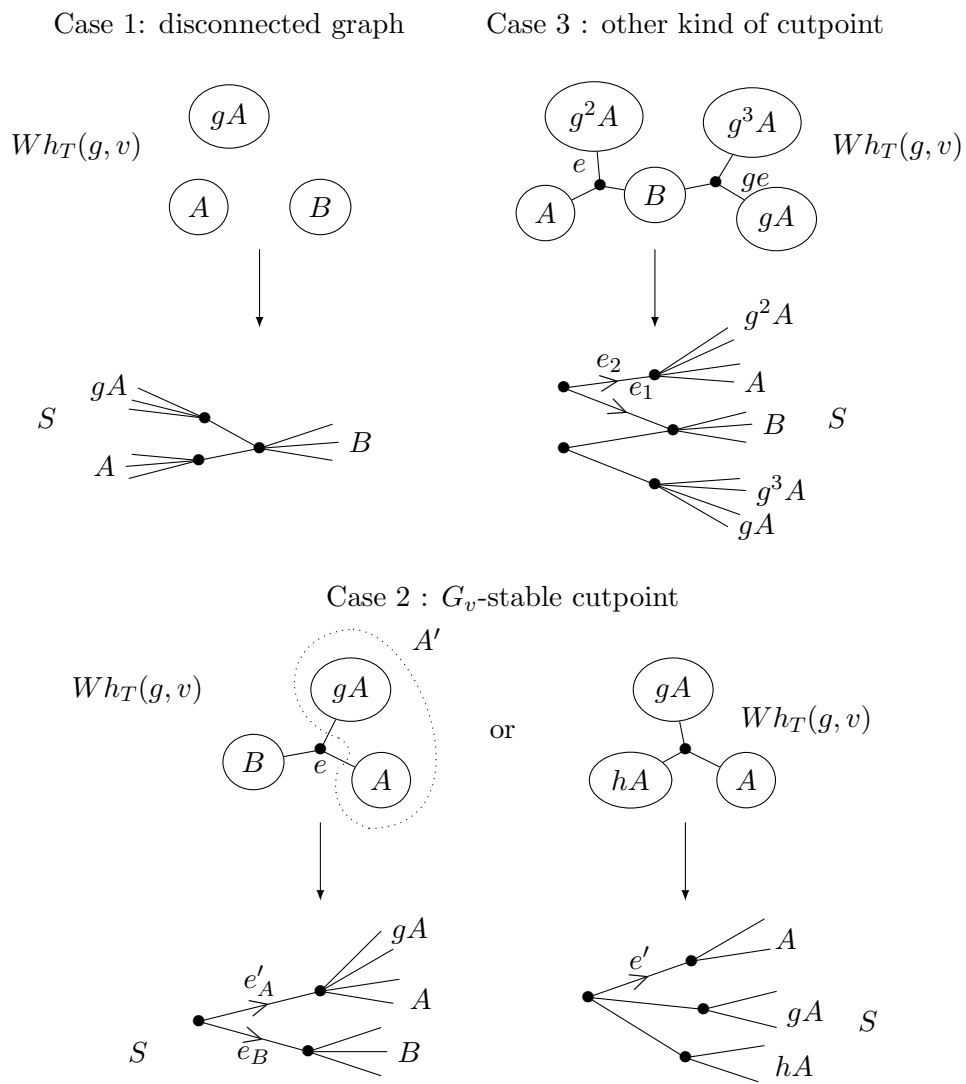

 \emph{Case 1 : }
  Suppose that the Whitehead graph $Wh_T(\G,v)$ is disconnected and that the stabilizer of a connected component is an allowed edge group of $\D$. In this case, denote the connected components by $C_1, \dots, C_n$ (note that two connected components might belong to the same orbit). Suppose that $\stab(C_1)$ is an allowed edge group.
 
  Since $C_1$ is a connected component of the Whitehead graph, we have $h\cdot C_1 \cap C_1 \neq \varnothing \Rightarrow h \in \stab(C_1)$. Let $S=T^{\stab(C_1), C_1}$ obtained by expansion, according to Lemma \ref{lem:constructionEclatement} with $f : S \to T$ the collapse map. We have $S \in \D^\A$ because $\stab C_1 \in \A$ and $T \in \D^\A$.
 
 Let $g \in \G$. As no translate of the axis of $g$ crosses a turn between any pair of distinct connected components $C_i$ and $C_j$, no translate of $\axe_S(g)$ can cross the added edges, so all turns of $\axe_T(g)$ lift to $S$. By lemma \ref{lem:relevementLocal} the collapse is isometric on the axes so $\|g\|_S=\|g\|_T$.

 \bigskip
 
 In other cases, suppose that the Whitehead graph has no admissible connected component but has a cut point (necessarily admissible, according to lemma \ref{lem:arbreDual}). In that case the orbit of the cut point under $G_v$ cannot be the whole graph. The reason is that the graph is finite and without simple loop, and thus cannot have only cut points as vertices.
 
 \bigskip
 
 \emph{Case 2: } 
  Suppose that the Whitehead graph $W$ has a $G_v$-invariant admissible cut point. Then $G_v$ fixes an edge in $T$. Call $e$ such a cut point. Denote by $A$ a connected component of $W \setminus \{e\}$: its stabilizer is a subgroup of $G_v$ so it lies in $\A$.
 
  If $G_v\cdot A \neq W \setminus \{e\}$ denote by $B$ the complement of $G_v \cdot A$ in $W \setminus \{e\}$. The part $B$ is stable under the action of $G_v$. The subset $A' := G_v \cdot A$ is stable as well. Neither $A'$ nor $B$ are empty. If $\{f, f'\}$ is a turn of a translate of the axis of some $g \in \G$, then $\{f,f'\}$ is included in $G_v \cdot A \cup \{e \}$ or in $B \cup \{e\}$.
 
  We define a new tree $S$ as follows (see Figure \ref{fig:casWhitehead}). First we expand $T$ at the vertex $v$ by unattaching edges of $A'$, attaching an edge $e_{A'}$ to $v$ and re-attaching the edges of $A'$ to the other end of $e_{A'}$, which gives the expanded tree $T_1:=T^{G_v,A'}$ (see lemma \ref{lem:constructionEclatement} for notations). Similarly we unattach edges of $B$ to re-attach them on a new edge $e_B$ with origin $v$, which gives the tree $T_2 := T_1^{G_v, B}$. The lemma guarantees that $T_2$ belongs to $\D^\A$ since the stabilizer of the new edge is in $\A$. Finally we collapse the edge $e$ of $T_2$, which is a collapsible edge since its stabilizer is $G_v$ and its ends are not in the same orbit. Let $S$ be the resulting tree. It belongs to $\D^\A$. Folding $e_{A'}$ with $e_B$ is a type A fold and yields the original tree $T$.

  Let $g$ be an element of $\G$. Let us prove $\|g\|_S=\|g\|_T$. The pre-image of an edge by the fold $S \to T$ is a star. Every turn in $\axe_T(g)$ lifts in $S$: the only turns which do not lift are those of the sort $h\cdot\{A',B\}$. Yet such turns are never crossed by the axis by assumption on the Whitehead graph. According to lemma \ref{lem:relevementLocal} $f$ is isometric in restriction to $\axe_S(g)$ so the translation length of $g$ is the same in $T$ as in $S$.

  \medskip

  On the contrary, if $G\cdot A = W \setminus \{e\}$ then as $A$ is not the only connected component there exists $h \neq 1$ such that $h\cdot A\cap A = \varnothing$. The stabilizer of $A$ is then a subgroup $H \subsetneq G_v$ and $H,A$ satisfy the conditions of the expansion lemma \ref{lem:constructionEclatement} with $H$ allowed as an edge stabilizer. We perform an expansion at the vertex $v$ as follows: for $u \in G_v/H$, 
  unattach $uA$ and re-attach it on a new edge $ue'$ with origin $v$.
  This yields the expanded tree $T^{H,A}$ which belongs to $\D^\A$. We then get $S$ by collapsing $e$. Since $e$ is collapsible we have $S \in \D^\A$. When one folds the edges $ue'$ with $u \in G_v/H$, one gets $T$; the fold is of type C.
 
  Again, for every $g \in \G$, all turns in $\axe_T(g)$ lift to $S$ so lemma \ref{lem:relevementLocal} guarantees $\|g\|_S=\|g\|_T$.

 \bigskip
 
 \emph{Case 3 : }
  Suppose that $W$ does not have any $G_v$-invariant cut point. Denote by $e$ a cut point and by $W_0$ the connected component of $W$ containing $e$. Let $A$ be a connected component of $W_0 \setminus \{e\}$ which does not contain any element of $G_v \cdot e$. Such a component exists by definition of an admissible cut point. By remark \ref{rem:stabilisateur} $A$ contains a vertex. Remember that its stabilizer is an allowed edge group (remark \ref{rem:stabilisateur}) and is a subgroup of $G_e$.
  
  Denote by $B$ the complement of $G_v \cdot A$ in $W \setminus G_v\cdot \{e\}$; $B$ may be empty and may intersect $W \setminus W_0$. 
  Again $B$ is stable by $G_v$. We also define $A':=G_e \cdot A$.
  
  Since $G_v \cdot e$ has at least two elements, $\{ e \} \cup A'$ is a proper subset of $\E_v$, so even when $B$ is empty, we may use lemma \ref{lem:constructionEclatement} to do the following expansions.
 
  See figure \ref{fig:depliage3} for a closer illustration of the case. First we do an expansion at vertex $v$: we partition the set of edges into $B \sqcup \bigsqcup_{h \in G_v /G_e} (\{e\} \cup A')$. We get the tree $T_1:=T^{G_e, \{e\} \cup A'}$ with notations of lemma \ref{lem:constructionEclatement}: we replace the vertex $v$ by a star with $|G_v/G_e|$ branches. Then we attach $B$ to the centre of the star, and edges in $h \cdot(\{e\} \cup A')$ to $h \cdot G_e$.

  We call $e_1$ the edge joining $\{e\} \cup A'$ to the centre of the star (which we will still call $v$). We call $w$ the origin of $e_1$.
  
  Then we perform a second expansion at $w$ which is the origin of $e_1$, $\bar e$ and of the edges of $A'$ and has stabilizer $G_e$. The tree $T_2 := T_1^{G_e, A'}$ may be described as follows: we unattach the edges of $A'$ then re-attach them on a new edge $e_2$ with origin $w$ (see figure \ref{fig:depliage3}).
 
  Finally consider the collapse $S = T_2 / \sim_e$. Since $G_e \in \A$, $T_1, T_2$ and $S$ are in $\D^\A$. There is an application $S \to T$ which sends $e_1$ and $e_2$ on $e$ and it is the (type B) fold of $e_1$ with $e_2$.

 \begin{figure}
  \centering
  %\usetikzlibrary{shapes}

\begin{tikzpicture}[scale=0.7]
%Wh schématique
\begin{scope}[xshift=4cm]
	\node[draw, ellipse] (A) at (-1,2) {$\; \vphantom g A' \;$};
	\node[draw, ellipse,text height=0.8cm, text width=0.5cm] (vide) at (0,0) {};
	\node[draw, ellipse] (gA) at (-1,-2) {$gA'$};
	
	\node[draw, ellipse,text height=1.2cm, text width=0.5cm] (reste) at (2,0) {};

	\coordinate (e) at (-0.1,1.34);
	\draw (e)  node {$\bullet$} node[left=0.15cm] {$e$} -- (vide);
	\draw (e) -- (A);

	\coordinate (ge) at (-0.1,-1.4);
	\draw (ge)  node {$\bullet$} node[left=0.1cm] {$ge$} -- (vide);
	\draw (ge) -- (gA);

	\draw[very thin, dotted] (-2,3) node[above right] {$W_0$} rectangle (1,-2.8);
	\draw[very thin, dotted] (-0.7,1.2) rectangle (2.7,-1.2) node [below left] {$B$};
	\node (leg) at (0.5,-3.5) {Shape of $W=Wh_T(g,v)$};
\end{scope}

%vue en haut au départ
\begin{scope}[yshift=-8.7cm]
	\coordinate (v) at (0,0);

	\draw (v) node [above] {$v$}  node {$\bullet$}-- +(-2.5,1) node [midway, sloped] {$>$} node[midway, above right] {$e$} node {$\bullet$} node [above] {$v'$};
	\draw (v) -- +(-2.5,-1) node [midway, sloped] {$>$} node[midway, below right] {$ge$} node {$\bullet$} node [above] {$gv'$};

	\draw (v) -- +(2.8,1.8) node[ above right] {$A'$};
	\draw (v) -- +(2.5,2.3);
	\draw (v) -- +(3.4,0.8);
	\draw (v) -- +(3.7,0.4) node [right] {$B$};
	\draw (v) -- +(3.6,-0.2);
	\draw (v) -- +(2.8,-1.8) node[ below right] {$gA'$};
	\draw (v) -- +(2.5,-2.3);

	\node[text width=5.5cm] (leg) at (0.5,-3.5) {Picture of $v$ in $T$};
\end{scope}

%premier éclatement

\begin{scope}[xshift=9.5cm, yshift=-8cm]
	\coordinate (v) at (0,0);
	\coordinate (w) at (-1,1);
	\coordinate (gw) at (-1,-1);

	\draw (w)  node {$\bullet$} node [above] {$w$}-- +(-2,1) node [midway, sloped] {$>$} node[midway, above right] {$e$} node {$\bullet$} node [above] {$v'$};
	\draw (gw) -- +(-2,-1) node [midway, sloped] {$>$} node[midway, below right] {$ge$} node {$\bullet$} node [above] {$gv'$};

	\draw (w) -- +(2.3,1.8) node[ above right] {$A'$};
	\draw (w) -- +(2,2.3);
	\draw (v)node [above] {$v$}  node {$\bullet$} -- +(2.4,0.8);
	\draw (v) -- +(2.7,0.4) node [right] {$B$};
	\draw (v) -- +(2.6,-0.2);
	\draw (gw)   node {$\bullet$} node [above left] {$gw$} -- +(2.3,-1.8) node[ below right] {$gA'$};
	\draw (gw) -- +(2,-2.3);

	\draw (v) -- (w) node [midway, sloped] {$>$} node[midway, above right] {$e_1$};
	\draw (v) -- (gw)node [midway, sloped] {$>$} node[midway, below right] {$ge_1$};

	\node[text width=5.2cm] (leg) at (0.5,-4.2) {$T_1$};
\end{scope}

% et le deuxième éclatement
\begin{scope}[xshift=0cm, yshift=-18cm]
	\coordinate (v) at (0,0);
	\coordinate (w) at (-1,1);
	\coordinate (gw) at (-1,-1);
	\coordinate (u) at (0,2);
	\coordinate (gu) at (0,-2);

	\draw (w)  node {$\bullet$} node [above] {$w$}-- +(-2,1) node [midway, sloped] {$>$} node[midway, above right] {$e$} node {$\bullet$} node [above] {$v'$};
	\draw (gw)  node {$\bullet$} node [above left] {$gw$}   -- +(-2,-1) node [midway, sloped] {$>$} node[midway, below right] {$ge$} node {$\bullet$} node [above] {$gv'$};

	\draw (u)  node {$\bullet$} -- +(1.3,0.8) node[ above right] {$A'$};
	\draw (u) -- +(1,1.3);
	\draw (v)node [above] {$v$}  node {$\bullet$} -- +(2.4,0.8);
	\draw (v) -- +(2.7,0.4) node [right] {$B$};
	\draw (v) -- +(2.6,-0.2);
	\draw (gu)   node {$\bullet$} -- +(1.3,-0.8) node[ below right] {$gA'$};
	\draw (gu) -- +(1,-1.3);

	\draw (v) -- (w) node [midway, sloped] {$>$} node[midway, above right] {$e_1$};
	\draw (v) -- (gw)node [midway, sloped] {$>$} node[midway, below right] {$ge_1$};

	\draw (u) -- (w) node [midway, sloped] {$>$} node[midway, above] {$e_2$};
	\draw (gu) -- (gw) node [midway, sloped] {$>$} node[midway,right] {$ge_2$};

	\node (leg) at (0.9,-4) {$T_2$};
\end{scope}

%terminé : résultat
\begin{scope}[xshift=9cm, yshift=-18cm]
	\coordinate (v) at (0,0);
	\coordinate (w) at (-2.5,1);
	\coordinate (gw) at (-2.5,-1);
	\coordinate (u) at (0,2);
	\coordinate (gu) at (0,-2);

	\draw (u)  node {$\bullet$} -- +(1.8,0.8) node[ above right] {$A'$};
	\draw (u) -- +(1.5,1.3);
	\draw (v)node [above] {$v$}  node {$\bullet$} -- +(3.4,0.8);
	\draw (v) -- +(3.7,0.4) node [right] {$B$};
	\draw (v) -- +(3.6,-0.2);
	\draw (gu)   node {$\bullet$} -- +(1.8,-0.8) node[ below right] {$gA'$};
	\draw (gu) -- +(1.5,-1.3);

	\draw (v) -- (w) node {$\bullet$}  node [midway, sloped] {$>$} node[midway, above right] {$e_1$};
	\draw (v) -- (gw)node {$\bullet$}  node [midway, sloped] {$>$} node[midway, below right] {$ge_1$};

	\draw (u) -- (w)  node {$\bullet$} node [above] {$v'$} node [midway, sloped] {$>$} node[midway, above] {$e_2$};
	\draw (gu) -- (gw) node {$\bullet$} node [above] {$gv'$} node [midway, sloped] {$>$} node[midway,right] {$ge_2$};

	\node (leg) at (0.5,-4) {$S$};
\end{scope}

%fléchage
%\draw[-latex, very thick] (4,0) -- (5,0);
\draw[latex-, very thick] (4.3,-9) -- (5.5,-9) node [midway, above]{Collapse of $e_1$}; %1
\draw[-latex, very thick] (4.3,-18) -- (5.5,-18)  node [midway, above]{Collapse of $e$};%3

%\draw [-latex, very thick] (10,-4) -- (10,-5);
\draw[latex-, very thick] (5.5, -13.5)--(4.3,-14.5)  node [midway, above, sloped]{Collapse of $e_2$};%2

\end{tikzpicture}
  \caption{Steps of the unfolding in the case 3, where the orbit of the cut point has several elements. Arrows represent collapses.} \label{fig:depliage3} 
 \end{figure}
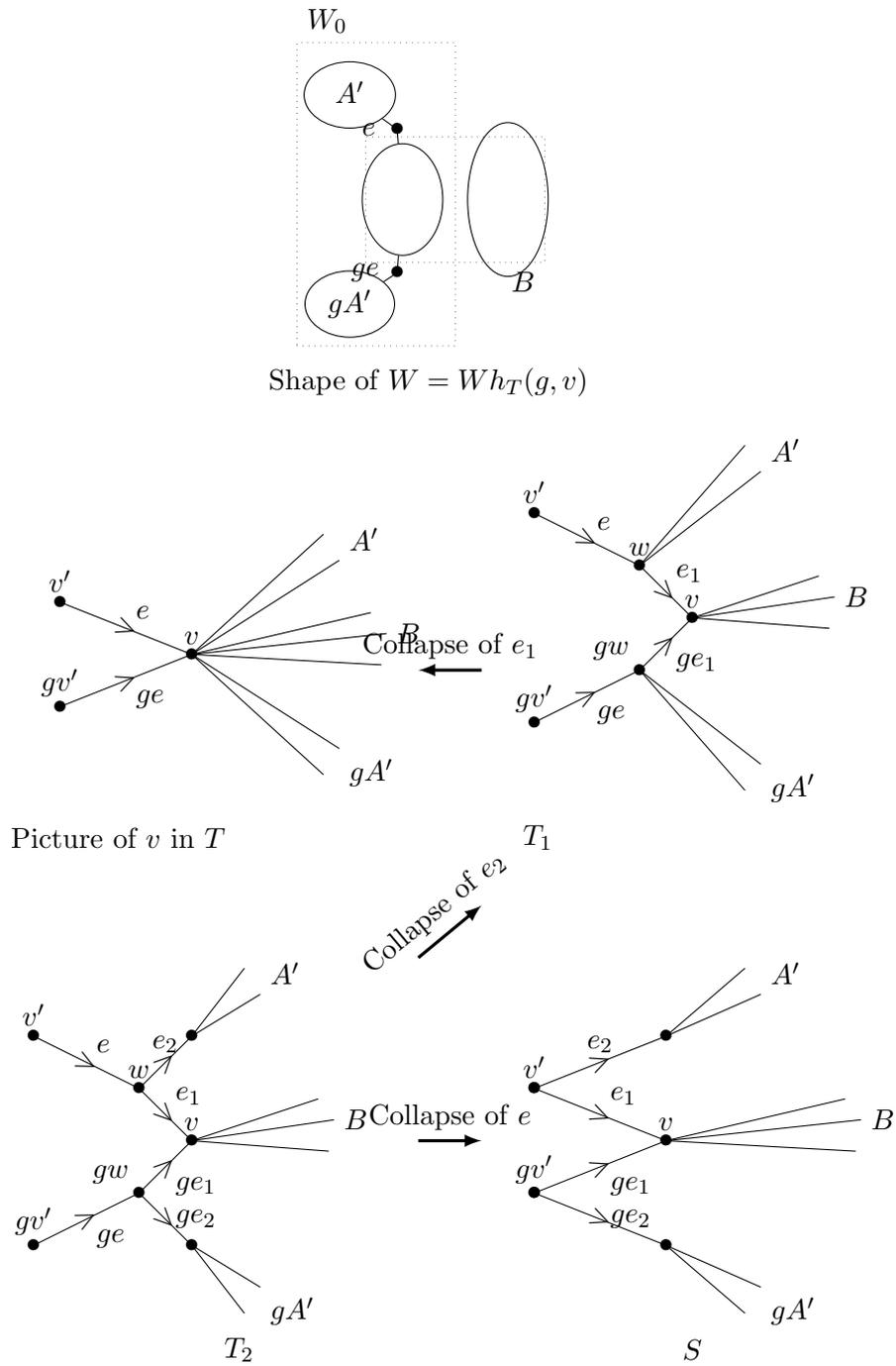

  Let $g$ be an element of $\G$.
  The only turns of $T$ at $v$ which may be crossed by translates of the axis of $g$ are those of the sort $h \cdot \{A,A\}$, $h \cdot \{A,e\}$, $ \{B,B\}$, $h \cdot \{B,e\}$, $h \cdot \{e,h'\cdot e\}$ with $h,h' \in G_v$. 
  All these turns lift to $S$. According to lemma \ref{lem:relevementLocal}, the whole axis lifts isometrically so $\|g\|_T = \|g\|_S$.
 
 \bigskip
 
  We have proved the existence of $S$ and $f$ in all cases where the Whitehead graph has an admissible cut.

 \bigskip
 
 Conversely, suppose there exists $S$ and $f: S \to T$ non-injective, sending edge to edge or edge to vertex, such that $\|g\|_S=\|g\|_T$ for every $g \in \G$. According to \cite{BestvinaFeighnBounding}, this application may be considered as a composition of collapses and folds.
 
 Let us consider only the last collapse or last fold. Since neither folds nor collapses can increase translation length, this application satisfies the second condition of the lemma. We may then suppose $f$ is either a fold or a collapse, which simplifies the proof. In both cases $f$ is $1$-Lipschitz. The assumption about translation lengths implies that for every $g \in \G$, $f$ is isometric on the axis of $g$. Therefore all turns in $\axe_T(g)$ lift to $S$.
 
 If $f$ is a collapse, as it does not change deformation space, it is a quasi-isometry, so connected components of the subforest collapsed by $f$ are bounded.
 
 Let $v \in T$ be such that the subtree $f^{-1}(v)$ is not reduced to a point. As $f$ is a collapse, the pre-image of any edge in $T$ is a single edge in $S$. In the Whitehead graph $Wh_T(\G,v)$, any two vertices joined by an edge correspond to edges of $T$ in the same connected component of $\overline{S\setminus f^{-1}(\{v\})}$. Otherwise $\axe_S(g)$ contains an edge collapsed by $f$. Therefore the Whitehead graph has at least as many connected components as $\overline{S\setminus f^{-1}(\{v\})}$ which is not connected.
 
 Let $e$ a collapsed edge in $S$ whose image is $v$. As $e$ is collapsible, it has an end $w$ such that $G_w = G_e$ and such that $w$ is terminal in $f^{-1}(v)$. The vertex $w$ belongs to the boundary of $S \setminus f^{-1}(v)$ since $S$ has no valence $1$ vertex. Every connected component of $S\setminus f^{-1}(\{v\})$ whose boundary is $w$ has stabilizer included in $G_e$. The stabilizers of all corresponding components in $Wh_T(\G,v)$ are then  subgroups of $G_e$, so they are allowed edge groups. Thus the Whitehead graph has an admissible cut.
 
 \bigskip
 
 If $f$ is a fold, it is defined by two edges of $S$ with same origin $w$: call them $e_1$ and $e_2$. Call their endpoints $v_1$ and $v_2$. Call $e'$ the edge of $T$ which is the image of $e_1$ and $e_2$, $w'$ its initial vertex and $v'$ its terminal vertex (which is the image of $v_1$ and $v_2$). The vertices $v'$ and $w' $ are in different orbits as we supposed that $T/G$ is without simple loop.

 We will prove that $Wh_T(\G,v') \setminus \{e'\}$ is disconnected and the stabilizer of at least one of its connected components $E_1$ is in $\A$. Lemma \ref{lem:arbreDual} states that this 
 implies that $Wh_T(\G,v')$ has an admissible cut.

\begin{figure}
  \centering
  \begin{tikzpicture}[scale=0.92]
%\draw [dotted] (0,10) rectangle (12.6,-10);

\begin{scope}[xshift=0.6cm]
\coordinate (w) at (0,0);
\coordinate (v1) at (1.8,1);
\coordinate (v2) at (2,0);
\coordinate (v3) at (1.8,-1);

\draw (w) node {$\bullet$} -- (v1) node {$\bullet$} node [midway, sloped] {$>$} node [midway, above] {$e$};
\draw (w) node {$\bullet$} -- (v2) node {$\bullet$} node [midway, sloped] {$>$} node [midway, above right] {$ae$};
\draw (w) node {$\bullet$} -- (v3) node {$\bullet$} node [midway, sloped] {$>$} node [midway, below] {$a^2e$};

\draw (v1) -- +(0.8,0.5);
\draw (v1) -- +(1,-0.1);
\draw (v2) -- +(1,0.3);
\draw (v2) -- +(1,-0.3);
\draw (v3) -- +(1,0.1);
\draw (v3) -- +(0.8,-0.5);

\draw (w) -- +(-0.6, 0.2);
\draw (w) -- +(-0.6, -0.2);

%ellipse
\draw[dotted] (2.7,1.2) ellipse [x radius=0.3,y radius=0.5] node [above right] {$A$};

% flèche vers le bas 
\draw [-latex, very thick] (1,-1.5) -- (1,-2.5) ;
\node[text width=2cm] (leg) at (2.5, -2) {Folding of $e$ with $ae$};

%le bas
\coordinate (ww) at (0,-4);
\coordinate (vv) at (2, -4);

\draw (ww) node {$\bullet$} -- (vv) node {$\bullet$} node [midway, sloped] {$>$} node [midway, above] {$e$};
%arêtes de droite
\draw (vv)  node [below left] {$v$}-- +(0.5,0.9);
\draw (vv) -- +(0.7,0.7);
\draw (vv) -- +(1,0.15);
\draw (vv) -- +(1,-0.15);
\draw (vv) -- +(0.8,-0.7);
\draw (vv) -- +(0.5,-0.9);

%arêtes de gauche
\draw (ww) -- +(-0.6, 0.2);
\draw (ww) -- +(-0.6, -0.2);

%ellipse basse
\draw[dotted] (2.6,-3.2) ellipse [x radius=0.25,y radius=0.35] node [above right] {$A$};

%les axes
%le vert
\draw [dashed,  green, thick] (-0.5,-3.7) .. controls (0,-3.95) .. (1,-3.95) .. controls (2,-3.95)  .. (2.4,-3) ;
\draw [dashed,  green, thick] (-0.5,0.3) .. controls (-0.1,0.1) and (0.1,0.1) .. (0.9,0.6) .. controls (1.7,1)  .. (2.6,1.6) ;

%le rouge
\draw [dashed, red, thick] (3,-0.4) .. controls (2,-0.06) .. (1,-0.05) .. controls (0,0) .. (1,-0.5) .. controls (1.8,-0.95) .. (2.8,-0.85);
\draw [dashed, red, thick] (3,-4.2) .. controls (2.1,-4.05) .. (2.9,-4.7);

% le graphe de Whitehead
\coordinate (p) at (0.5,-7.5);
\draw (p) node {$\bullet$} -- +(1.7,1) coordinate (p1) node {$\bullet$};
\draw (p) -- +(1.9,0.7) node {$\bullet$} coordinate (p2);
\draw  (p1)-- (p2);

\draw (p)  -- +(2.1,0.3) coordinate (p1) node {$\bullet$};
\draw (p) -- +(2.1,-0.2) coordinate (p2) node {$\bullet$};
\draw  (p1)-- (p2);

\draw (p) -- +(1.7,-1) coordinate (p1) node {$\bullet$};
\draw (p)  -- +(1.9,-0.7) coordinate (p2) node {$\bullet$};
\draw  (p1)-- (p2);

\node (leg) at (1.6,-9.3) { Shape of $Wh_T(g,v)$};
\node (leg) at (1.6,2.3) { First case};
\end{scope}

\begin{scope}[xshift=5.5cm]
%\node (leg) at (1.6,-9.3) { Shape of $Wh_T(g,v)$};
\node (leg) at (2.8,2.3) { Second case, with condition 1 or 2};

%dessin du haut
\coordinate (w) at (0,0);
\coordinate (v1) at (1.8,0.7);
\coordinate (v2) at (1.8,-0.7);
\coordinate (v3) at (-0.1,-2);
\coordinate (v4) at (0.5,-2);

\draw (w) node {$\bullet$} -- (v1) node {$\bullet$} node [midway, sloped] {$>$} node [midway, above] {$e_1$};
\draw  (w)  -- (v2) node {$\bullet$} node [midway, sloped] {$>$} node [midway, below] {$e_2$};
\draw[dotted] (w) -- (v3) node {$\bullet$} node [midway, sloped] {$>$} node [midway, below left] {$ae_2$};
\draw[dotted] (w)  -- (v4) node {$\bullet$} node [midway, sloped] {$>$} node [midway, below right] {$ae_1$};

\draw (v1) -- +(0.6,0.5);
\draw (v1) -- +(0.8,-0.1);
\draw (v2) -- +(0.8,0.1);
\draw (v2) -- +(0.6,-0.5);
\draw (v2) -- +(0.78,-0.2);

\draw [dotted](v3) -- +(0.2,-0.5);
\draw[dotted] (v3) -- +(-0.2,-0.5);
\draw [dotted](v4) -- +(0.2,-0.5);
\draw [dotted](v4) -- +(-0.2,-0.5);

\draw (w) -- +(-0.6, 0.2);
\draw (w) -- +(-0.6, -0.2);

%ellipse
\draw[dotted] (2.6,1) ellipse [x radius=0.3,y radius=0.5] node [above right] {$A$};
\draw[dotted] (2.5,-0.9) ellipse [x radius=0.3,y radius=0.5] node [below right] {$B$};

%le bas
\coordinate (ww) at (0,-4);
\coordinate (vv) at (2, -4);
\coordinate (vv3) at (0.5,-5.5);

\draw (ww) node {$\bullet$} -- (vv) node {$\bullet$} node [midway, sloped] {$>$} node [midway, above] {$e$};
\draw[dotted] (ww) -- (vv3) node {$\bullet$} node [midway, sloped] {$>$} node [midway, below left] {$ae$};

%arêtes de droite
\draw (vv)  node [below left] {$v$}-- +(0.5,0.9);
\draw (vv) -- +(0.7,0.5);

\draw (vv) -- +(0.8,-0.5);
\draw (vv) -- +(0.5,-0.9);
\draw (vv) -- +(0.7,-0.75);

\draw [dotted](vv3) -- +(0.2,-0.5);
\draw[dotted] (vv3) -- +(-0.2,-0.5);
\draw [dotted](vv3) -- +(0.1,-0.6);
\draw[dotted] (vv3) -- +(-0.1,-0.6);
\draw[dotted] (vv3) -- +(0,-0.6);

%arêtes de gauche
\draw (ww) -- +(-0.6, 0.2);
\draw (ww) -- +(-0.6, -0.2);

%ellipse basse
\draw[dotted] (2.6,-3.3) ellipse [x radius=0.25,y radius=0.35] node [above right] {$A$};
\draw[dotted] (2.65,-4.7) ellipse [x radius=0.25,y radius=0.35] node [below right] {$B$};

% le graphe de Whitehead
\coordinate (p) at (0.5,-7.5);
\draw (p) node {$\bullet$} -- +(1.7,1) coordinate (p1) node {$\bullet$} node [right] {$A$};
\draw (p) -- +(1.9,0.7) node {$\bullet$} coordinate (p2);
\draw  (p1)-- (p2);

\draw (p) -- +(1.7,-1) coordinate (p1) node {$\bullet$};
\draw (p)  -- +(2,-0.7) coordinate (p2) node {$\bullet$} node [right] {$B$};
\draw (p)  -- +(2,-0.3) coordinate (p3) node {$\bullet$};
\draw  (p1)-- (p2);
\draw  (p1)-- (p3);
\draw  (p3)-- (p2);

\end{scope}
%dernier cas

\begin{scope}[xshift=11cm]
	\coordinate (v) at (0,0);
	\coordinate (w) at (-1.5,0.4);
	\coordinate (gw) at (-1.5,-0.4);
	\coordinate (u) at (0,1);
	\coordinate (gu) at (0,-1);

%moustaches de droite
	\draw (u)  node {$\bullet$}  node [above ] {$v'$}-- +(0.8,0.3);% node[ above=-0.1cm] {$A$};
	\draw (u) -- +(0.5,0.6);
	\draw (v)node [above] {$v$}  node {$\bullet$} -- +(1,0.5);% node[ below=0.6cm] {$B$};
	\draw (v) -- +(1.1,0.1)  ;
	\draw (v) -- +(1,-0.2);
	\draw (gu)   node {$\bullet$} -- +(0.8,-0.3);% node[ below =0.1cm] {$aA$};
	\draw (gu) -- +(0.5,-0.6);

	\draw (v) -- (w) node {$\bullet$}  node [midway, sloped] {$>$} node[midway, above right=-0.1cm] {$e_2$};
	\draw (v) -- (gw)node {$\bullet$}  node [midway, sloped] {$>$} node[midway, below right=-0.05cm] {$ae_2$};

	\draw (w)  node {$\bullet$} node [above] {$w$} -- (u) node [midway, sloped] {$>$} node[midway, above=0cm] {$e_1$} ;

	\draw (gw) node {$\bullet$} node [above] {$aw$} -- (gu) node [midway, sloped] {$>$} node[midway,below=0.1cm] {$ae_1$};

	 \draw (w) -- +(-0.5, 0.2);
	 \draw (w) -- +(-0.5, -0.2);
	 \draw (gw) -- +(-0.5, 0.2);
	 \draw (gw) -- +(-0.5, -0.2);

%ellipse
\draw[dotted] (0.6,1.5) ellipse [x radius=0.3,y radius=0.5] node [above right] {$A$};
\draw[dotted] (1,0.15) ellipse [x radius=0.3,y radius=0.5] node [below right] {$B$};
\draw[dotted] (0.6,-1.5) ellipse [x radius=0.3,y radius=0.5] node [below right] {$aA$};

%dessin du bas 
	\coordinate (v) at (-0.2,-4);
	\coordinate (w) at (-1.3,-3.5);
	\coordinate (gw) at (-1.3,-4.5);

	\draw (v)  -- +(1,1) ;
	\draw (v) -- +(1.1,0.8);

	%E2
	\draw (v) node [above= 0.2cm] {$v$}  node {$\bullet$} -- +(1.2,0.3);
	\draw (v) -- +(1.4,0.1) ;
	\draw (v) -- +(1.3,-0.1);
	%aE1
	\draw (v)   -- +(1,-0.7);
	\draw (v) -- +(0.8,-0.9);

	\draw (v) -- (w) node {$\bullet$}  node [midway, sloped] {$>$} node[midway, above right] {$e$};
	\draw (v) -- (gw)node {$\bullet$}  node [midway, sloped] {$>$} node[midway, below right] {$ae$};

	 \draw (w) -- +(-0.5, 0.2);
	 \draw (w) -- +(-0.5, -0.2);
	 \draw (gw) -- +(-0.5, 0.2);
	 \draw (gw) -- +(-0.5, -0.2);

%ellipse
\draw[dotted] (0.7,-3.2) ellipse [x radius=0.3,y radius=0.3] node [above right] {$A$};
\draw[dotted] (1,-3.9) ellipse [x radius=0.3,y radius=0.4] node [below right] {$B$};
\draw[dotted] (0.6,-4.7) ellipse [x radius=0.3,y radius=0.3] node [below right] {$aA$};

% le graphe de Whitehead
\coordinate (p) at (-1.2,-7.2);
\coordinate (pp) at (-1.2,-7.8);
\draw (p) node {$\bullet$} -- +(1.7,1) coordinate (p1) node {$\bullet$} node [right] {$A$};
\draw (p) -- +(1.9,0.7) node {$\bullet$} coordinate (p2);
\draw  (p1)-- (p2);

\draw (pp) node {$\bullet$} -- +(1.7,-1) coordinate (p1) node {$\bullet$} node [right] {$A$};
\draw (pp) -- +(1.9,-0.7) node {$\bullet$} coordinate (p2);
\draw  (p1)-- (p2);

\draw (p) -- +(2,-0.7) coordinate (p1) node {$\bullet$};
\draw (p)  -- +(1.8,-0.3) coordinate (p2) node {$\bullet$} node [right=0.2cm] {$B$};
\draw (p)  -- +(2,0.1) coordinate (p3) node {$\bullet$};
\draw  (p1)-- (p2);
\draw  (p1)-- (p3);
\draw  (p3)-- (p2);

\draw (pp) -- (p1);
\draw (pp) -- (p2);
\draw (pp) -- (p3);

\end{scope}
\end{tikzpicture}
  \caption{Summary of all different cases of folds which may occur and corresponding shape of associated Whitehead graph}  \label{fig:casDepliages}
 \end{figure}
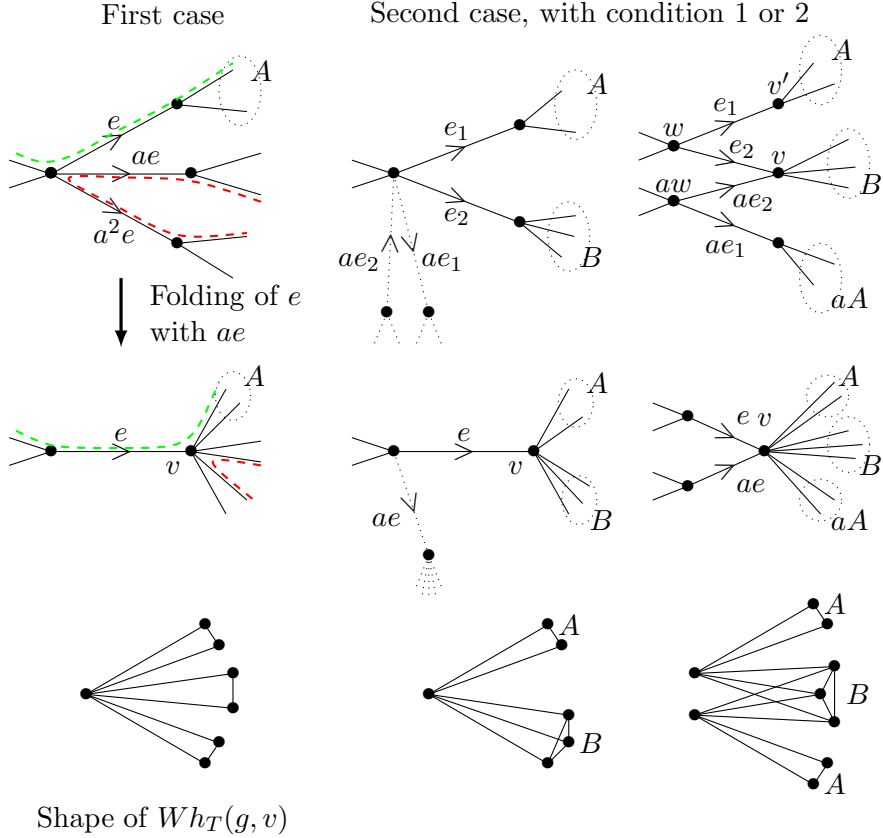

 Figure \ref{fig:casDepliages} recaps all different cases of folds and associated shapes of graphs.
 
 \begin{figure}[h]
  \centering
  \begin{tikzpicture}[scale=1]
\draw (0,0) node {$\bullet$} node [above] {$w$} -- (2,1) node [midway, sloped, above] {$e_1$} node {$\bullet$} node [above] {$v_1$};
\draw (0,0)-- (2,-1) node [midway, sloped, above] {$e_2$} node {$\bullet$} node [above] {$v_2$};

 \draw[-latex ] (3,0) -- (4,0);

\draw (5,0)node {$\bullet$} node [above] {$w'$}  -- (7.25,0) node [midway, above] {$e'$} node {$\bullet$} node [above] {$v'$};

%aretes accessoires
\draw [dotted] (0,0) -- (-0.9,0.2);
\draw [dotted] (0,0) -- (-0.9,-0.2);
\draw [dotted] (0,0) -- (-0.9,0.5);
\draw [dotted] (0,0) -- (-0.9,-0.5);

\draw[dotted] (2,1) -- (2.9,1.2);
\draw[dotted] (2,1) -- (2.9,0.8);
\draw[dotted] (2,1) -- (3,1);
\draw[dotted] (2,-1) -- (2.9,-1.15);
\draw[dotted] (2,-1) -- (2.9,-0.85);

\begin{scope}[xshift=5cm]
\draw [dotted] (0,0) -- (-0.9,0.2);
\draw [dotted] (0,0) -- (-0.9,-0.2);
\draw [dotted] (0,0) -- (-0.9,0.5);
\draw [dotted] (0,0) -- (-0.9,-0.5);

\draw [dotted] (2.3,0)--(3,0);
\draw [dotted] (2.3,0)--(3,0.2);
\draw [dotted] (2.3,0)--(3,-0.2);
\draw [dotted] (2.3,0)--(3,0.4);
\draw [dotted] (2.3,0)--(3,-0.4);
\end{scope}

\end{tikzpicture}
  
  \caption{A fold as described.}
 \end{figure}
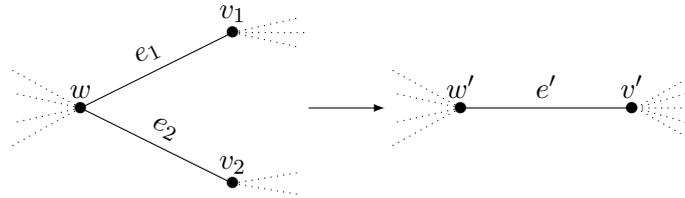
 
 Three different kind of folds may occur, which correspond to cases A, B and C.
  \begin{enumerate}
  \item $G_{e_i}= G_{v_i}$ for $i\in \{1,2\}$ and $e_1, e_2$ lie in different orbits (type A)
  \item $G_{v_1} \subset G_{e_2}$, up to permutation of indices, and $e_1, e_2$ lie in different orbits (type B)
  \item $e_2 = h e_1$ for some $h \in G_v$ and $G_{e_i} = G_{v_i}$ (type C)
 \end{enumerate}

In the three cases, define $\tilde E_1$ to be the set of edges of $S$ with origin $v_1$, except $\bar e_1$. 
Let $E_1$ be the image of $\tilde E_1$ in $T$.
Define $E_2$ similarly.

 Since all turns represented in $Wh_T(\G,v')$ lift to $S$, any edge of $Wh_T(\G,v')$ with one endpoint in $E_1$ joins $E_1$ to itself or to $\bar e'$. In particular, $E_1$ and $E_2$ are in distinct connected components of $Wh_T(\G,v') \setminus \{e'\}$.
 Therefore $Wh_T(\G,v')$ has an admissible cut.
\end{proof}

\subsection{Proof of the theorem} \label{subsec:proofOfTheo}

\begin{lem} \label{lem:eclatementArbre}
 (Expansion of non-cyclic vertex groups)
 Let $R \in \bar \D^\A$ whose certain vertices $v_1, \dots, v_k$ (in different orbits) have stabilizers $H_1, \dots, H_k$ some special factors. Let $T \in \D$ and $T_i$ the minimal subtree for $H_i$ in $T$. There exists a $G$-tree $S \in \D^\A$ and a map $f : S \to T$ and a collapse $\pi : S \to R$ such that:
 \begin{itemize}
  \item the image by $\pi$ of the collapsed subforest is $G \cdot \{v_1, \dots, v_k\}$,
  \item  for all $i \in \{1, \dots, k\}$, $f_{|\pi^{-1}(v_i)}$ is an isomorphism to $T_i$.
 \end{itemize}
\end{lem}

In other words, it is possible to blow up $R$ by replacing $v_i$ by $T_i$.

A proof of this result is given in \cite[Proposition 2.2]{GuirardelLevitt17}. The key assumption is the fact that all edge groups in $R$ are elliptic in $T$ so we can attach edges of $R$ to the subtrees which replace the vertices. Here we suppose that $R$ is the result of the collapse of some $\tilde R \in \D$, so its edge groups are also edge groups in $\tilde R$ and are elliptic in any tree in $\D^\A$.

\begin{proof}[Proof of theorem \ref{theo:Whitehead}.]

Let $G$ be a GBS group and let $\G:= \{g_1, \dots, g_k\}$ be a finite collection of loxodromic elements of $G$. Let $T \in \D^\A$. Suppose that $\G$ is simple with respect to $\D^\A$, that is, there exists a non-trivial $G$-tree $R \in \bar \D^\A$ such that every $g_i \in \G$ fixes a vertex  $v_i \in R$.

For every $1 \leq i \leq k$, let $T_i$ be the minimal $G_{v_i}$-invariant subtree of $T$. We obtain a new tree $S\in \D^\A$ by applying lemma \ref{lem:eclatementArbre} starting from $R$ so that for every $i \in \{1, \dots, k\}$, the vertex $v_i$ is replaced by a copy of $T_i$.

Let $f: S \to T$ be the map given by Lemma \ref{lem:grapheWhEtDepliage}. If $f$ is injective then it is an isomorphism (surjectivity is obtained by minimality of the image for the action of $G$). In that case, for every $i \in \{1, \dots, k\}$, the axis of $g_i$ avoids images of edges in $S$ coming from edges in $R$. Therefore some Whitehead graph at an end of such an edge has an isolated vertex and has an admissible cut.

When $f$ is not injective, we showed Lemma \ref{lem:grapheWhEtDepliage} that there is a vertex at which the Whitehead graph has an admissible cut, which proves the theorem.
\end{proof}

\section{Algorithm} \label{sec:algo}

This part is dedicated to a theorem which states that one can decide algorithmically whether a collection of loxodromic elements of $G$ is simple with respect to some $\D^\A$, or not.

Before stating the theorem, let us explain how we deal with the set $\A$ of allowed edge groups. Let $\Gamma$ be a graph of groups representing $G$; it has finitely many vertices $v_1, \dots, v_n$. 

\begin{defi}
 For $H$ a subgroup of $G$ we denote by $\A _H$ the family of subgroups $\{A \cap H, A \in \A \}$. Since $\A$ is stable by taking subgroups, $\A_H$ is a subfamily of $\A$.
 
 For each $v \in V(\Gamma)$ let $I_v$ be a family of positive integers. We say that $(I_v)_{v\in V(\Gamma)}$ \emph{represents} $\A$ if for every $v \in V$, $\A_{|G_v}$ is the set of all subgroups of $G_v$ whose index is a multiple of an element in $I_v$.
\end{defi}

\begin{ex}
 
 \begin{itemize}
  \item $I_v = \{1\}$ for each $v \in \Gamma$ represents $\A_\Z$.
  \item Suppose $\Gamma$ is reduced. Let $I_v$ be the set of absolute values of labels at $v$, then $(I_v)_{v \in \Gamma}$ represents $\Amin$. 
  \item  Given $\Gamma$ a graph of groups for $G$, we can choose $\A$ to be the set of bi-elliptic subgroups in $\Gamma$. If so, $\A$ is represented by $(I_v)_{v \in V(\Gamma)}$ where $I_v$ is the set of all labels at $v$.
 \end{itemize}
\end{ex}
In the sequel we assume that $\A$ is represented by a family $(I_v)_{v \in V(\Gamma)}$ of finite sets.

Here is the theorem which we will prove.
\begin{theo} \label{theo:algoWhitehead}
 There is an algorithm which takes as input:
 \begin{itemize}
  \item a graph of cyclic groups $\Gamma_G$ representing $G$
  \item a marking $G \to \pi_1(\Gamma_G,v_0)$
  \item a finite family of finite subsets $(I_v)_{v \in V(\Gamma_G)}$ representing a family $\mathcal A$ of subgroups of $G$
  \item a finite collection of loxodromic elements $\G \subset G$.
 \end{itemize}
  which decides whether there exists a system of proper special factors $\mathcal H$  of $G$ with respect to $\D^\A$ such that $\G \preceq \mathcal H$, and then returns such a system when it exists.
\end{theo}

We show that the construction of Whitehead graphs is algorithmic, and that the construction of the map of Lemma \ref{lem:grapheWhEtDepliage} is algothmic too. Finally, after giving a description of the algorithm, we prove that it terminates.

\subsection{Paths in graphs of groups, fundamental group}

In this section we give some general definitions and results about GBS groups. For algorithmic purposes it is more convenient to work with graphs of groups than with trees. 

Let $G$ be a GBS group. We suppose it is given by a graph of cyclic groups $\Gamma_G$. From this graph we deduce a presentation $\langle s_1, \dots, s_n| r_1, \dots, r_m \rangle$. We will use this presentation to define a marking in other graphs of groups below.

Let $\Gamma$ be any graph of cyclic groups. It is given as a set of vertices and edges, a label $\lambda(e)$ for each edge $e$ and a set of generators $(a_v)_{v \in V(\Gamma)}$ for all vertex groups.

We define $\pi_1(\Gamma)$ as a subgroup of the Bass group $B(\Gamma)$ like in \cite[Definition (a), part 5.1]{SerreArbresAmalgames}. The group $B(\Gamma)$ has the following presentation:
generators are 
\begin{itemize}
 \item elements $a_v$ indexed by vertices of $\Gamma$, where $a_v$ is to be thought of as a generator of the vertex group $G_v \simeq \Z$
 \item elements $t_e$ indexed by oriented edges in $\Gamma$.
\end{itemize}

Relations are the following:
\begin{enumerate}
 \item for all edge $e \in \Gamma$ we have $t_e= t_{\bar e }^{-1}$
 \item for all $e \in \Gamma$ with initial vertex $u$ and terminal vertex $v$, with labels $\lambda(e)=p$ and $\lambda(\bar e)=q$, we have $a_{v}^{q}= t_e a_{u}^{p} t_e^{-1}$
\end{enumerate}

\bigskip

A \emph{path $\alpha$ in the graph of groups} is a pair $\alpha=(w,\gamma)$ where 
\begin{itemize}
 \item $\gamma$ is a path $v_0, e_1, \dots, e_n, v_n$ in the underlying graph of $\Gamma$, where the $v_i$ are vertices and $e_i$ are edges such that $t(e_i)= v_i=o(e_{i+1})$,
 \item $w$ is a word $a_0 t_{e_1} a_1 \dots t_{e_n} a_n$ where $a_i \in G_{v_i}$ for all $i \in \{0, \dots, n\}$.
\end{itemize}

We denote by $[\alpha]$ the element of $B(\Gamma)$ represented by $w$.

The \emph{length} of a path $(w,\gamma)$ is the number of edges in $\gamma$. Its initial and terminal vertices are $v_0$ and $v_n$ respectively.

Let $\alpha =(a_0 t_{1} \dots a_n; v_0, e_1, \dots ,v_n)$ and $\alpha' =(a'_0 t'_1 \dots a_m'; v'_0, e'_1, \dots, v'_m)$ be two paths in $\Gamma$ such that $v_n=v'_0$. The \emph{concatenation} of $\alpha$ and $\alpha'$ is the path 
\[
 \alpha \cdot \alpha'= (a_0 t_1 \dots t_n b t'_1 a'_1 \dots a'_m; v_0, e_1, \dots, v_n, e'_1, \dots, v'_m)
\] 
where $b$ is the element of $G_{v_n}$ equal to $a_n a'_0$. 

Let $v \in \Gamma$. A \emph{loop based at $v$} in $\Gamma$ is a path with initial and terminal vertices equal to $v$.

To a loop $\alpha = (w, \gamma)$, one associates the corresponding element $[\alpha]:= [w] \in B(\Gamma)$.

If $\alpha, \alpha'$ are loops in $\Gamma$ based at $v$ then $[\alpha \cdot \alpha'] = [\alpha]\cdot [\alpha']$.

Fix a vertex $v$ in $\Gamma$. The fundamental group $\pi_1(\Gamma,v)$ is the subgroup of $B(\Gamma)$ consisting of the elements of $B(\Gamma)$ associated to loops based in $v$. 

\begin{rema} \label{rem:mot}
 If a word in the generators of $B(\Gamma)$ corresponds to a path in the graph, then this loop is unique. We really mean the word as a sequence fo letters and not the corresponding element of $B(\Gamma)$. Thus the word is sufficient to describe a path in the graph of groups, and we will use the word on its own when the description in terms of edges and vertices is not needed.
\end{rema}

A \emph{marking} of $\Gamma$ is a map $\{s_i, 1\leq i \leq n\} \to \pi_1(\Gamma,v)$ which associates a loop based at $v$ in $\Gamma$ to each generator of $G$, such that it induces an isomorphism $G \to \pi_1(\Gamma,v)$.

\bigskip

For every $g \in \G$, given the expression of $g$ as a word in the generators $\{s_i, 1 \leq i \leq n\}$, we can determine a loop in $\Gamma$ based at $v$ which represents $g$.

We call a path in $\Gamma$ given by a word $w$ \emph{reduced} if no subword in $w$ is of the form $t_e a t_{\bar e}$ with $a \in i_{e}(G_e)$. Note that here vertex groups are cyclic so determining if $a$ belongs to $i_{e}(G_e)$ boils down to a question of divisibility.

We call a loop $\alpha$ \emph{cyclically reduced} if the concatenation $\alpha \cdot \alpha$ is reduced.

\bigskip

We may modify the marking by the following process. Given another vertex $v'$ and a path $\alpha$ in the graph of groups $\Gamma$ from $v$ to $v'$, there is an isomorphism $\pi_1(\Gamma,v) \simeq \pi_1(\Gamma,v') $ defined by 

\begin{align*}
  \sigma_\alpha : &\pi_1(\Gamma,v)  \to \pi_1(\Gamma,v') \\
 & [h] \mapsto [\bar \alpha h  \alpha]
\end{align*}

\begin{lem} \label{lem:reduction}
 A path $h=(w,\gamma)$ can be reduced algorithmically, i.e. there is an algorithm which finds a reduced path $h'$ such that $[h] = [h']$ in $B(\Gamma)$.
 
 A loop can be cyclically reduced algorithmically: for any loop $\alpha$, one can find a cyclically reduced loop $\alpha'$ and a path $\beta$ such that $[\alpha] = [\bar \beta \alpha' \beta]$.
\end{lem}

The proof is standard and straightforward. We leave it to the reader.

\begin{lem} \label{lem:reduitLongueurNulle}
 Suppose that $\gamma$ is a reduced loop based at $v$ and that $[\gamma] \in G_v$. Then $\gamma$ has length $0$. 
\end{lem}
The proof for this fact is in \cite[5.2, Theorem 11]{SerreArbresAmalgames}.

With the elements above we define the universal cover $T_{\Gamma,v}$ of the graph of groups $\Gamma$. It is a graph defined as follows. The set of vertices is
\[ 
 \tilde V = \{ \text{paths in } \Gamma \text{ with initial vertex } v \}/ \sim
\]

where $\gamma \sim \gamma'$ if $\gamma$ and $\gamma'$ have the same terminal vertex $v_i \in \Gamma$ and $[\gamma]^{-1} [\gamma'] \in G_{v_i}$. Denote by $[\gamma]_V$ the vertex associated to the path $\gamma$.
The group $\pi_1(\Gamma,v)$ acts on $\tilde V$ by left concatenation.

Note that checking whether two paths define the same vertex boils down to checking whether a path in $\Gamma$ can be reduced to a length zero path, by lemma \ref{lem:reduitLongueurNulle}, so it is algorithmic.

Let $\tilde v = [1]_V$. It is a lift for the base point.

The oriented edges of $(T,v)$ are defined as follows: 
\[ 
 \tilde E =  \left \{ (\alpha,a t_e)/ \begin{array}{l}
                                 \alpha \text{ path in } \Gamma \text{ from } v \text{ to } v', e \in \E_{v'}, a \in G_{v'}
                               \end{array}
  \right \}/ \sim
\]
where $\E _{v'}$ is the set of edges with origin $v'$. The equivalence relation $\sim$ is defined by $(\alpha,a t_e) \sim (\alpha',a' t_{e'})$ if and only if 
$e = e'$ and  $ a^{-1} \alpha^{-1} \alpha' a' \in i_{\bar e}(G_e)$. 
The origin of this edge is $[\alpha]_V$ and its terminus is $[\alpha \cdot a t_e]_V$.

Denote by $[(\alpha, a t_e)]_E$ the equivalence class of $(\alpha, a t_e)$. The group $\pi_1(\Gamma, v)$ acts on $\tilde E$ by left concatenation of $\alpha$. We have the relation $[(\alpha, a t_e)]_E = [(\alpha \cdot a, t_e)]_E$ for all $a \in G_{v'}$. The edge with opposite orientation is $\overline{[(\alpha, a t_e)]_E} = [\alpha \cdot at_e, t_{\bar e}]$.

\medskip

The graph defined is a tree.
The quotient of $(T,v)$ under the action of $\pi_1 (\Gamma, v)$ is $\Gamma$; the projection of $[\gamma]_V$ is the last vertex of $\gamma$ (see \cite{SerreArbresAmalgames} for a proof).

Note that for different base points $v,v'$ the universal covers $T_{\Gamma,v}$ and $T_{\Gamma, v'}$ are isomorphic. Given $\alpha$ joining $v$ to $v'$, we define an isomorphism between $T_{\Gamma,v}$ and $T_{\Gamma, v'}$ by $[\gamma]_V \mapsto [\bar \alpha \cdot \gamma]_V$ (this is well-defined since it does not depend on the representative $\gamma$ chosen).
A marking $G \to \pi_1(\Gamma,*)$ yields an action of $G$ on $T_{\Gamma, *}$.
If we identify $\pi_1(\Gamma,v)$ with $\pi_1(\Gamma,v')$ with the isomorphism $\sigma_\alpha$ as above, the isomorphism between the trees is $G$-equivariant.

\bigskip

\begin{lem} \label{lem:reducedGeodesic}
 Let $\alpha$ be a path in $\Gamma$ with first vertex $v$. It can be lifted to a path in $T$ with first vertex $\tilde v$ and last vertex $p_\alpha$, the equivalence class of $\alpha$ in $\tilde V$.
 
 Suppose $\alpha$ is reduced. Then the distance between $\tilde v$ and $p_\alpha$ in $(T, v)$ is equal to the length of $\alpha$.
\end{lem}

\begin{proof}
 We follow the proof of \cite[Proposition 2.6]{BeekerMultiple}.
 
 First of all we construct the path between $\tilde v$ and $p_\alpha$.
 Denote by $n$ the length of $\alpha$. Write (as a word) $\alpha =a_0 t_{e_1} a_1 \dots a_n$, and  $\alpha_i :=a_0 t_{e_1} a_1 \dots a_i$ for $i \leq n$. 
 
 Let $\tilde v_i = [\alpha_i]_V$ and $\tilde e_i = [(\alpha_{i-1}, t_{e_{i}})]_E$.
 
 The path $\tilde v, \tilde e_1, \tilde v_1, \dots , \tilde e_n, \tilde v_n = g \cdot \tilde v$ is a path in $T$. Let us show that it is a geodesic, that is, no consecutive edge are opposite.
 
 Suppose $\overline{ \tilde e_{i+1}} = \tilde e_{i}$ for some $i \in \{0, \dots, n\}$. Then \[\overline{[(\alpha_{i}, t_{e_{i+1}})]_E} = [(\alpha_i \cdot t_{e_{i+1}}, \bar t_{e_{i+1}})]_E = [(\alpha_{i-1} , t_{e_{i}})]_E\]
 so $\bar e_i = e_{i+1}$ and $t_{e_i} a_i \bar t_{e_i}  \in i_{e_i}(G_{e_i})$. Since the word is reduced this leads to a contradiction. Thus the lift of the path $\alpha$ is a geodesic in $T$.
\end{proof}

\begin{lem} \label{lem:domFonda}
 Let $\Gamma$ be a labelled graph and $g$ be an element in $\pi_1(\Gamma,v)$ represented by a cyclically reduced loop $\gamma$. Let $T$ be the universal cover of $\Gamma$ defined as above.
 
 If the loop has length zero, then $g$ is elliptic. Otherwise $g$ is loxodromic and the path in $T$ defined by $\gamma$, joining $\tilde v$ to $g \tilde v$, is a fundamental domain of the axis of $g$ in $T$.
 
 In particular, for each loxodromic element $g \in G$, one can compute a fundamental domain of its axis.
\end{lem}

\begin{proof}
 Let $\gamma$ be a cyclically reduced loop in $\Gamma$ such that $[\gamma] = g$. Let $\tilde v$ be the lift of the base point $v$ in $T$. Since $\gamma$ is cyclically reduced, lemma \ref{lem:reducedGeodesic} ensures that $[\tilde v,  g^2 \cdot \tilde v]$ is a geodesic. Therefore $\tilde v$ belongs to the axis of $g$ in $T$. As a result the path $[\tilde v, g \cdot \tilde v]$ is a fundamental domain of the axis of $g$.
\end{proof}

\subsection{Algorithmicity of Whitehead graph computation and unfoldings}

\begin{lem} \label{lem:computeLink}
 Let $\Gamma$ be a marked graph of groups for $G$. Let $T$ be its universal cover.
 Let $x$ be a vertex in $T$ represented as $x = [\gamma_x]_V$ for some path $\gamma$ in $\Gamma$. Let $G_x$ be the stabilizer of $x$, not to be confused with vertex groups in $\Gamma$. The subgroup $G_x$ is generated by $a_x := \gamma_x a_{\pi(x)} \bar \gamma_x$ and one can compute
 \begin{itemize}
  \item the link $lk(x)$
  \item the action of $a_x$ on $lk(x)$.
 \end{itemize}
\end{lem}
\begin{proof}
 Since $T$ is a tree, $lk(x)$ is the collection of edges of $T$ with initial vertex $x$. Denote by $\pi$ the quotient map $T \to \Gamma$. 
 Any edge of $T$ with origin $x$ has a unique representative of the form $(\gamma_x, a_{\pi(x)}^k t_e)$ with $e \in lk(\pi(x))$ and $0 \leq k < \lambda(e)$.

 All such edges can be listed algorithmically since the indices of edge groups in vertex groups are all finite.
 
 \medskip
 
 For every $y:= [\gamma_x, a_{\pi(x)}^k t_e]_E \in lk(x)$ we have $a_x \cdot y = [\gamma_x,  a_{\pi(x)}^{k+1 \text{ mod } \lambda(e)} t_e]_E $.
\end{proof}

Let $g \in G$ be a loxodromic element. With the input of a graph of groups and a loop for $g$ we can compute the Whitehead graphs $Wh_T(\{g\},v)$:

\begin{lem} \label{lem:computeWh}
 Let $\Gamma$ be a graph of cyclic groups. Let $(w, \gamma)$ be a cyclically reduced loop in $\Gamma$ based in $v$ representing some $g \in \pi_1(\Gamma)$. Let $T$ be the universal cover of $\Gamma$  at basepoint $v$ and let $x$ be a vertex in $T$. The computation of the Whitehead graphs $Wh_T(\{g\},x)$ is algorithmic. 
\end{lem}
\begin{proof}
 By lemma \ref{lem:computeLink} we can compute the link of $x$ and the action of $G_x$ on it.

 By lemma \ref{lem:domFonda} we may compute a fundamental domain of the axis of $g^2$ (or equivalently a pair of consecutive fundamental domains of $g$). 
 All orbits of turns crossed by the axis of $g$ appear in this segment.
 For every turn $\tau$ of $\axe(g)$ based at a point in the orbit of $x$ we can find a pair of edges in $lk (x)$ forming a turn in the same orbit, and using the action of $G_x$, we can find all turns at $x$ in the orbit of $\tau$.
 
 The elements of $lk(x)$ form the vertices of $Wh_T(g,x)$ and the turns computed above are edges.
\end{proof}

\begin{cor}
 Let $\G \subset G$ be a finite collection of loxodromic elements of $G$. Then the computation of the Whitehead graphs $Wh_T(\G,v)$ is algorithmic. 
 \qed
\end{cor}

\begin{lem}
 Given a Whitehead graph $Wh_T(\G,v)$, the action of $G_v$ on it and the set $I_v$, one can decide algorithmically whether it has an admissible cut.
\end{lem}
\begin{proof}
 Finding connected components in a finite graph is algorithmic. The set $I_v$ allows to check the admissibility of connected components by calculating stabilizers.
\end{proof}

\bigskip

Lemma \ref{lem:grapheWhEtDepliage} states that given a Whitehead graph with an admissible cut there exists a non-injective map $f:S \to T$ which preserves translation length of every $g\in \G$. In the proof we actually gave a construction of such a map. We are going to prove that this construction can be done algorithmically. First we show that collapses and expansions can be done algorithmically:

%===à suivre===%

%===fin de la partie algorithmique===%

\begin{lem} \label{lem:algoEcrase}
 There is an algorithm which takes as input a graph of groups $\Gamma$, a collapsible edge $\epsilon=vw$ of $\Gamma$ such that $G_{w}= G_\epsilon$, a marking $\psi: G \to \pi_1(\Gamma,v)$
 and outputs
 \begin{itemize}
  \item a marked graph of groups $(\Gamma',\psi')$ whose universal cover is the tree $T'$ obtained by the collapse of the orbit $\epsilon$ of the universal cover $T$ of $(\Gamma, \psi)$ 
  \item an isomorphism $\phi: \pi_1(\Gamma,v) \to \pi_1(\Gamma',v')$ such that $\phi \circ \psi= \psi'$. 
 \end{itemize}
\end{lem}
\begin{proof}
 First construct the graph of groups $\Gamma'$: we collapse the collapsible edge $\epsilon$ in $\Gamma$. Since $G_{w}= G_\epsilon$ we have $\lambda(\bar \epsilon) = \pm 1$. We create the new graph of group $\Gamma'$ by deleting $\epsilon$ and $w$ and redefining any edge with origin $w$ by attaching it to $v$ instead. The label of such edges is multiplied by $\pm \lambda(\epsilon)$. There is a map between the underlying graphs $f :\Gamma \to \Gamma'$ which sends each vertex to the corresponding vertex in $\Gamma'$ (which we will write with a $'$) and sends $w$ to $v'$, and sends any edge except $\epsilon$ to the corresponding edge and sends $e$ to the vertex $v'$.
 
 Now let us build a morphism $\phi:B(\Gamma)\to B(\Gamma')$:
   \begin{align*}
   a_u & \mapsto a_{u'}  \text{ for } u \neq w \\
   a_{w} & \mapsto a_{v'}^{\pm\lambda(\epsilon)} \\
   t_e & \mapsto t_{e'} \text{ for } \epsilon \neq e \\
   t_\epsilon & \mapsto 1  
   \end{align*}
 This morphism sends any path in the graph of groups $\Gamma$ to a path in $\Gamma'$. In particular it induces a morphism $\pi_1(\Gamma, v) \to \pi_1(\Gamma',v')$. Moreover the induced morphism is an isomorphism between both fundamental groups. We define $\psi'=\phi \circ \psi$. It is a marking on $\Gamma'$.
 
 The morphism $\phi$ induces a map $(T,v) \to (T',v')$ because it sends paths to paths and preserves the equivalence relations defining trees. 
 This map is a collapse map and it is $G$-equivariant for the markings $\psi$ and $\psi'$.
 
The graph $\Gamma'$ and the marking $\psi'$ can be computed from $\Gamma$ and $\psi$, which proves the lemma.
\end{proof}

We have a similar result for expansions. It is slightly more difficult since while collapses may be defined from data in the quotient, the definition of an expansion requires some information in the tree.

\begin{lem} \label{lem:algoEclate}
  There is an algorithm which takes as input
  \begin{itemize}
   \item a marked graph of groups $(\Gamma, \psi)$
   \item a vertex $v \in \Gamma$
   \item a subset $S \subsetneq E(\tilde v)$, where $\tilde v$ is the basepoint of the universal cover $(T,v)$, such that $\forall g \in G_{\tilde v}, \quad g S \cap S \neq \varnothing \Rightarrow g \in \stab(S)$
  \end{itemize}
   and outputs
 \begin{itemize}
  \item a marked graph $(\Gamma',\psi')$ whose universal cover  is $T'=T^{\stab(S),S}$ (see lemma \ref{lem:constructionEclatement})
  \item an isomorphism $\phi: \pi_1(\Gamma,v) \to \pi_1(\Gamma',v')$ such that $\phi \circ \psi= \psi'$. 
 \end{itemize}
\end{lem}

\begin{proof}
 Like in the proof of the collapse in lemma \ref{lem:algoEclate}, we first construct an oriented graph, then construct a map between $B(\Gamma)$ and $B(\Gamma')$ which induces a map between trees.
 
 Let $\Gamma'$ be the labelled graph obtained as follows. The vertices of $\Gamma'$ are the vertices of $\Gamma$ along with another vertex $w'$. To distinguish them from the vertices of $\Gamma$ we write their names with a $'$. The edges are redefined according to the following rule (see figure \ref{fig:redef_aretes}). Let $n:= [G_{\tilde v}: \stab(S)]$. Let $\pi: T \to \Gamma$ be the quotient map.
 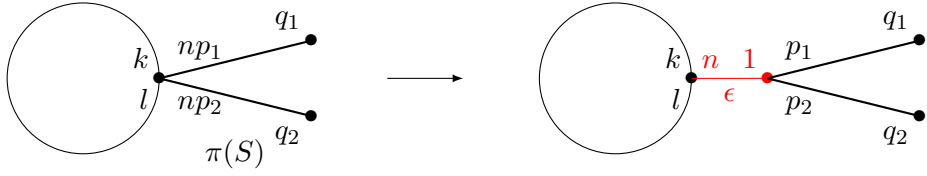
\begin{figure}
  \centering
  \begin{tikzpicture}
 \draw (0,0) node {$\bullet$} arc (0:360:1) node [above left] {$k$} node [below left] {$l$};
 \draw[thick] (0,0) node [above right=0.1cm] {$np_1$} -- (2,0.5) node [above left] {$q_1$} node {$\bullet$};
 \draw[thick] (0,0) node [below right=0.1cm] {$np_2$} -- (2,-0.5) node [below left] {$q_2$} node {$\bullet$};
 
 \draw (1, -1) node {$\pi(S)$};
 
 \draw[-latex] (3,0) -- (4,0);
 
 \begin{scope}[xshift = 7cm]
  
 \draw (0,0) node {$\bullet$} arc (0:360:1) node [above left] {$k$} node [below left] {$l$};
 \draw[color = red] (0,0) node [above right] {$n$} -- (1,0) node [above left] {$1$} node [midway, below] {$\epsilon$} node {$\bullet$};
 
 \draw[thick] (1,0) node [above right=0.1cm] {$p_1$} -- (3,0.5) node [above left] {$q_1$}node {$\bullet$};
 \draw[thick] (1,0) node [below right=0.1cm] {$p_2$} -- (3,-0.5) node [below left] {$q_2$} node {$\bullet$};
 \end{scope}

\end{tikzpicture}
  \caption{Redefinition of edges} \label{fig:redef_aretes}
 \end{figure}

 \begin{enumerate}
  \item Add an edge $\epsilon$ with origin $v'$ and terminal vertex $w'$ with label $n$ near $v'$ and $1$ near $w'$. 
  \item For $e \in E(\Gamma)$ not in $\pi(S)$, keep $e$ with same origin and same label $\lambda(e)$.
  \item For $e \in E(\Gamma)$ in $\pi(S)$ we have $o(e)=v$. Redefine its origin to be $w'$ and its label to $\lambda(e)/n$.
 \end{enumerate}
 This defines the origin of all oriented edges of $\Gamma'$. We then define $t(e)= o(\bar e)$ for all $e \in E(\Gamma')$.
 
 Let $f: \Gamma' \to \Gamma$ be the natural collapse map.
 
 For each edge $e$ with origin $v$ in $\Gamma$, let $\tilde e :=[(1,t_e)]_E$ be its standard lift in $T$ with origin $\tilde v$. If $e \in \pi(S)$ 
 choose $c(e)$ such that $e_S := a_v^{c(e)} \tilde e \in S$. Note that $\pi ^{-1}(e) \cap S = \langle a_v^n \rangle \tilde e_S$.

 Let $\chi: B(\Gamma') \to B(\Gamma)$ be the following morphism:
   \begin{align*}
   a_{u'} & \mapsto a_{u} \text{ for } u' \neq w'\\
   a_{w'} & \mapsto a_v^n\\
   t_{e'} & \mapsto t_{e} \text{ for } e \notin \pi(S) \cup \pi(\bar S) \\
   t_{e'} & \mapsto  a_v^{c(e)} t_{e} \text{ for }  e \in \pi(S)\\
   t_{\epsilon} & \mapsto 1 
   \end{align*}
   and define the morphism $\phi : B(\Gamma) \to B(\Gamma')$:
   \begin{align*}
   a_{u} & \mapsto a_{u'}\\
   t_e & \mapsto t_{e'} \text{ for } e \notin \pi(S) \cup \pi(\bar S) \\
   t_{e} & \mapsto  a_{v}^{-c(e)} t_\epsilon t_{e'}  \text{ for }  e \in \pi(S)
   \end{align*}
We can check that $\chi \circ \phi = \id _{B(\Gamma)}$. Both these morphism send paths to paths so they induce morphisms between fundamental groups. One also easily checks that $\phi \circ \chi_{|\pi_1(\Gamma',v')} = \id _{\pi_1(\Gamma',v')}$. Thus $\chi$ induces an isomorphism $\pi_1(\Gamma', v') \to \pi_1(\Gamma,v)$ and the restriction of $\phi$ is its inverse.

\medskip

Define the marking $\psi'= \phi \circ \psi$ of $\Gamma'$.

The morphism $\chi$ sends paths to paths, so it induces a map $\tilde f: (T', v') \to (T,v)$ between universal covers. This map is $G$-equivariant for the markings $\psi$ and $\psi'$.

\medskip

We want to prove that $\tilde f$ is the collapse map of the orbit $\epsilon$ and that $T'= T^{\stab S,S}$.

Let $\tilde v'= [1_{v'}]_V$ be the basepoint in $T'$. Let $\tilde w'= [1_{v'} t_{\epsilon} 1_{w'}]_V$ and $\tilde \epsilon = [(1_{v'}, t_\epsilon)]_E$.
In order to prove that $T'=T^{\stab(S),S}$ we need to show that the link of the vertex $\tilde w'$ is the pre-image of $S$ by $\tilde f$, together with the edge $\Bar{\tilde \epsilon}$.

Note that $\tilde f$ collapses the orbit $\epsilon$. To prove that it is a collapse map, we also need the following fact: $\tilde f$ is injective on the set of edges not in the orbit of $\epsilon$. Indeed, for $e \in E(T') \setminus G\cdot \epsilon$, $\phi \circ \chi (G_e) = G_e$. Thus two distinct edges in the same orbit cannot be sent to the same edge.

Note that edges in different orbits are also sent to images in different orbits. This proves that $\tilde f$ is injective on the set of edges. Thus it is injective on the interior of $T' \setminus G\cdot \tilde \epsilon$ and it is the collapse map associated to the orbit $\epsilon$.

Let $S'=lk(\tilde w) \setminus \{\tilde \epsilon\}$. Let us show that it is contained in the preimage of $S$. 
Let $\tilde e':=[(1_{v'} t_\epsilon a_{w'}^k, t_{e'})]_E$ be an edge in $S'$. By definition of the edges in $\Gamma'$, the corresponding orbit of edge $e \in \Gamma$ is in $\pi(S) \subset \Gamma$. Thus there is an edge $[(a_{w'}^{c(e)}, t_e)]_E \in S \subset T$. In view of the definition of $\chi$, the image of $\tilde e'$ by $\tilde f$ is $[(a_v^{nk},a_v^{c(e)} t_e)]_E=[(a_v^{nk+c(e)},t_e)]_E$. Since $a_v^{nk} \in \stab(S)$ we have $\tilde f (\tilde e) =[a_v^{nk}(a_v^{c(e)},t_e)]_E = a_v^{nk} \tilde e_S$ so it belongs to $S$.

\medskip

The number of edges in $S'$ is equal to the number of edges in $S$, and both can be computed with $[G_{\tilde v} : \stab S]$ and the labels. Since $\tilde f$ is injective on the set of edges and sends $S'$ to $S$, it induces a bijection between both sets so $\tilde f ^{-1} (S)= S'$. This proves that $T'$ is the tree $T^{\stab S, S}$.
\end{proof}

\begin{cor}\label{cor:lemmesEclateEcrase}
 Given $T$, $H$ and $S$ such as defined in lemma \ref{lem:constructionEclatement}, the tree $T^{H,S}$ can be constructed algorithmically.
 \qed
\end{cor}

\begin{lem} \label{lem:technique}
 There is an algorithm which takes as input
 \begin{itemize}
  \item a marked graph of groups $(\Gamma_T, \psi)$ with $\pi_1(\Gamma_T)\simeq G$, and universal cover $T$
  \item finite sets $I^T_v$ associated to each vertex $v \in V(\Gamma_T)$ representing a family of allowed edge groups $\mathcal A$ such that $T \in \D^{\mathcal A}$
  \item a Whitehead graph $Wh_T(\G,v)$ and an admissible cut of this graph for $\mathcal A$, where $\G$ is a finite collection of loxodromic elements of $G$
 \end{itemize}
 and gives
 \begin{itemize}
  \item a marked graph of groups $\Gamma_S$ and universal cover $S \in \D^{\mathcal A}$
  \item sets $I^S_v$ associated to vertices of $\Gamma_S$ representing the same collection of allowed subgroups $\mathcal A$
 \end{itemize}
 such that there exists a non-injective map $S \to T$, sending edge to edge or vertex, such that for every $g \in \G$, $\|g\|_S = \|g\|_T$. Moreover the map $S \to T$ can be chosen to be either a fold or a collapse. 
\end{lem}

\begin{proof}
 The proof of lemma \ref{lem:grapheWhEtDepliage} gives a construction of such a tree $S$ (or equivalently the associated graph of groups) by performing expansions and collapses on $T$. Lemmas \ref{lem:algoEcrase} and \ref{lem:algoEclate} compute a new marked graph of groups for a collapse or an expansion, along with an isomorphism between fundamental groups which is compatible with the markings.
 
 Finally we need to check that the sets which describe $\A$ can be computed, i.e. for any vertex group in $\Gamma_S$ we need to find its maximal subgroups which belong to $\A$. For any vertex group $H$ of $\Gamma_S$ the change of markings given by lemmas \ref{lem:algoEcrase} and \ref{lem:algoEclate} enables one to compute the image of $H$ as a vertex subgroup in $\Gamma_T$, so using the family $(I^T_v)_{v \in \Gamma_T}$ one can compute the family $(I^S_v)_{v \in \Gamma_S}$.
\end{proof}

\subsection{Description and termination of the algorithm}

 In this subsection we prove Theorem \ref{theo:algoWhitehead}. Using the criterion given by theorem \ref{theo:Whitehead} we give an algorithm deciding, given $G$ and $\G \subset G$ a finite collection of loxodromic elements, whether $\G$ is simple, and if so, determining a system of special factors $\mathcal H$ such that $\G \preceq \mathcal H$.

Note the following fact:

\begin{lem} \label{lem:sommetIsoleWh}
 If $Wh_T(\G,v)$ has an isolated vertex, then $\G$ is simple with respect to $\D$.
\end{lem}
\begin{proof}
 In the quotient $T/G$ the axis of $g$ avoids some edge which defines a subforest of $T/G$ containing the image of the axis of every $g \in \G$.
\end{proof}

Here is the description of the algorithm. We start with the tree $T_0 :=T_G \in \D$ corresponding to the graph of groups $\Gamma_G$ defining $G$. We may check immediately whether $G$ is a solvable Baumslag-Solitar group, in which case no proper special factor exists. We suppose it 
is not the case. Start with $n = 0$:
\begin{enumerate}
 \item Compute the axis of every $g \in \G$ in $T_n$. If some edge orbit does not intersect any axis, then its complementary subgraph is a subforest of $\Gamma$. It may have some components with elliptic fundamental group. Such components do not contain any element of $\G$. The set of components of the subforest with non-elliptic fundamental group gives a system of proper special factors $\mathcal H$ such that $\G \preceq \mathcal H$. The algorithm returns YES.
 \item Compute Whitehead graphs for all vertices of $T_n/G$ and check whether at least one of them has an admissible cut, using lemma \ref{lem:computeWh}. If none is found, then $\G$ is not simple and the algorithm returns NO.
 \item If we find $v \in T_n/G$ such that $Wh_{T_n}(\G,v)$ has an admissible cut, we compute a tree $T_{n+1}$ such that there is a $G$-equivariant $f: T_{n+1} \to T_n$ sending edge to edge or to vertex with $\|g\|_{T_{n+1}}=\|g\|_{T_n}$ for every $g \in \G$. Lemma \ref{lem:technique} ensures that this can be done. Start again step 1 with $T_{n+1}$.
\end{enumerate}

\begin{lem}
 The algorithm described above terminates. When it does, either it finds $T^*$ such that $\axe_{T^*} (g)$ does not cross some orbit of edges for every $g \in \G$, or it finds a proof that $\G$ is not simple.
\end{lem}

\begin{proof}
 At each iteration of the second step we replace $T_n$ by a tree $T_{n+1}$ obtained by an expansion or an unfolding. Expansions and unfoldings of type A or B increase the number of orbits of edges by one while unfoldings of type C do not change this number. Suppose by contradiction that the algorithm does not terminate, yielding an infinite sequence $(T_n)_{n \in \N}$. Lemma \ref{lem:pasdexplosion} implies that the number of edges in $T_n/G$ tends to infinity. Thus there exists $N \in \N$ such that the number of edges of $T_N/G$ is strictly greater than $\displaystyle \sum_{g \in \G} \|g\|_{T_N} = \displaystyle \sum_{g \in \G}\|g\|_{T_0}$. 
 Then there must be at least one orbit of edge avoided by the axis of every $g \in \G$, so the algorithm should have stopped at the $N$-th iteration: this is the contradiction we needed.
\end{proof}

\section{Decreasing sequences of special factors} \label{sec:suitesDecroissantes}

In this part, $G$ is a non-elementary GBS group. Let $\A$ be a family of allowed edge subgroups of $G$.
Let $H$ be a special factor with respect to $\D^{\A}$. The induced deformation space $\D^{\A}_{|H}$ is the deformation space of $H$-trees of which elliptic groups and allowed edge groups are those of $\D$ which are contained in $H$.

In particular $\D_{|H}$ is the cyclic deformation space for $H$. If $\A _{|H} := \{ A \cap H / A \in \A\}$, we have $\D^{\A}_{|H} = (\D_{|H})^{\A _{|H}}$.

\begin{lem}
 A subgroup $K \subset H$ is a special factor of $G$ with respect to $\D^{\mathcal A}$ if and only if $K$ is a special factor of $H$ with respect to $\D^{\mathcal A}_{|H}$.
\end{lem}
\begin{proof}
 If $T$ is a tree in $\bar \D^{\mathcal A}$ then its minimal $H$-invariant subtree $T_H$ is a tree in $\bar \D^{\mathcal A}_{|H}$. Conversely, any $H$-tree in $\D_{|H}$ is a subtree of some $G$-tree in $\D$.
 
 Suppose $K \subset H$ is a special factor of $G$, then it is a vertex stabilizer of some vertex $v$ in some $T \in \bar \D^{\mathcal A}$. The vertex $v$ must belong to the minimal $H$-invariant subtree $T_H \in \bar \D^{\mathcal A}_{|H}$ since $K$ is not an allowed edge group.
 Therefore $K$ is a special factor for $H$ with respect to $\D^{\mathcal A}_{|H}$.
 
 Now suppose $K$ is a special factor of $H$ with respect to $\D^{\mathcal A}_{|H}$: there is a collapse $T \to \bar T$ with $T \in \D^{\mathcal A}_{|H}$ and such that $K$ is a vertex stabilizer in $\bar T$. There is a tree $S \in \D^{\mathcal A}$ such that $T$ is a subtree of $S$ and there is a collapse $S \to \bar S$ such that the restriction to $T$ is the collapse $T \to \bar T$ and the collapse is an isometry on $S \setminus T$. Edge stabilizers belong to $\mathcal A$. Thus $K$ is a vertex stabilizer in $\bar S$ so it is a special factor for $G$.
\end{proof}

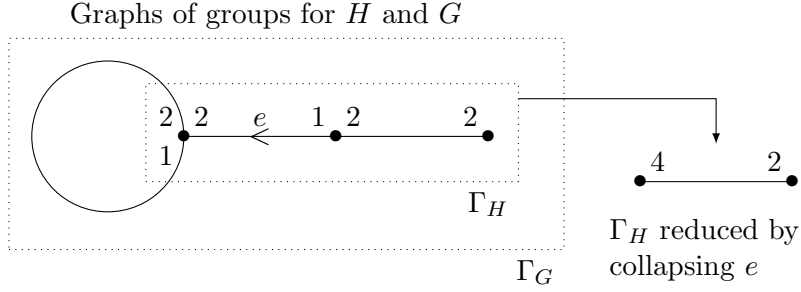
\begin{figure}
 \centering
 \begin{tikzpicture}[]
	%\draw[dotted] (-2.5,-10) rectangle (9.5,10);
	\draw (0,0) arc (0:360:1) node (x) {$\bullet$};
	\draw (x) node [above left]{$2$} node [below left] {$1$} node [above right] {$2$};
	\draw (0,0) --(2,0) node (y) {$\bullet$} node [midway] {$<$} node [midway, above ] {$e$} -- (4,0) node (z) {$\bullet$};
	\draw (y) node [above left]{$1$}  node [above right] {$2$};
	\draw (z) node [above left]{$2$};

%grand rectangle 
	\draw [dotted] (-2.3, 1.3) rectangle (5,-1.5) node [below left] {$\Gamma_G$};
	\draw [dotted] (-0.5, 0.7) rectangle (4.4,-0.6) node [below left] {$\Gamma_H$};
	\draw [-latex] (4.4,0.5) -| (7,-0.1);
	\draw (6,-0.6)  node (xx) {$\bullet$} -- (8,-0.6) node (zz) {$\bullet$};
	\draw (xx)  node [above right] {$4$};
	\draw (zz) node [above left]{$2$};
\node[text width=2.8cm] (leg) at (7,-1.5) {$\Gamma_H$ reduced by collapsing $e$};
\node[text width=7cm] (leggauche) at (2,1.6) {Graphs of groups  for $H$ and $G$};
\end{tikzpicture}
 \caption{Example: $H$ is a special factor of $G$ for $\D_{\text{red}}$ because $G_e$ fixes the edge of the loop after reduction. On the contrary $(\D_{\text{red}})_{|H}$ is not the reduced space for $H$ because $G_e$ does not fix any edge in reduced $H$-trees (for example the universal cover of the graph of groups on the right).} \label{fig:exempleDred}
\end{figure}

\begin{rema} \label{rem:facteursReduits}
 If $\mathcal A = \Amin$ the family of subgroups which are bi-elliptic in reduced trees of $\D$, we have $\D^{\mathcal A} = \D_{\text{red}}$. However, the elements of $\D^{\Amin}_{|H}$ are not necessarily reduced as $H$-trees of $\D_{|H}$. In fact some bi-elliptic groups appearing in some reduced $G$-trees might not be bi-elliptic in reduced $H$-trees (see figure \ref{fig:exempleDred} for an example). This means that $\Amin_{|H}$ \emph{is not} the family of subgroups of $H$ which are bi-elliptic in reduced $H$-trees of $\D_{|H}$.
\end{rema}

Given a family $\G$ of loxodromic elements of $G$, the algorithm of Theorem \ref{theo:algoWhitehead} that we described can be iterated in order to find a $\preceq$-decreasing sequence $(\mathcal H _n)$ of systems of special factors with respect to $\D^{\A}$ such that for all $n \in \N$, $\G \preceq \mathcal H_n$. 

One can ask whether there exists a minimal system of special factors $\mathcal H _{\min} $ such that $\G$ is $\mathcal H _{\min} $-peripheral, and if the iteration of the algorithm eventually finds such a system. The answer is yes. We introduce a complexity on special factors which enables us to prove the existence of a minimal factor and that the algorithm stops.

\begin{rema} 
 Additional operations used to iterate the Whitehead algorithm are themselves algorithmic. We need two things. The first one is to be able to take a (possibly non-connected) subgraph $\Gamma' \subset \Gamma$, compute the corresponding subgroups $G'_1, \dots, G'_k$ and the minimal subtrees of these subgroup. Indeed the connected components of $\Gamma'$ may have valence $1$ vertices, in which case they do not represent the minimal subtree for the corresponding $G'_i$. The second one is to deduce the replacements for $\D^\A$ in every factor of the system. It consists in keeping all allowed edge groups which are contained in the special factor $G'_i$ for every $i \in \{1, \dots, k\}$. In practice, for a factor $G'_i$ corresponding to a connected component $\Gamma'_i$, this means we keep every $I_v$ for which $v$ belongs to $\Gamma'_i$.
\end{rema}

For any non-elementary GBS group $H$ we denote by $b_1(H)$ the first Betti number of any graph of cyclic groups with fundamental group $H$. This is an invariant of trees in the cyclic deformation space for $H$. We denote by $M(H)$ the set of conjugacy classes of big vertex stabilizers with respect to $\Amin$. Define $m(H):=\#M(H)$. Vertex stabilizers are the same in all trees in $\D_{|H}$ so $m(H)$ is also well-defined.

We also introduce the following integer:
 \[
  \sigma(H) = \sum_{K \in M(H)} i(K)
 \]
where $i(K)$ is an integer which we define as follows and which is linked to the peripheral structure of $K$ (see \cite[Definition 4.10]{GuirardelLevitt07})

Define $i(K) = [K : K']$ where $K' = \langle g \in K | g \text{ bi-elliptic in some } T \in \D^{\Amin} \rangle$. This definition does not depend on any graph of groups for $H$. However one can compute $i(K)$ easily using the labels of a reduced graph of groups for $H$. Given such a graph, there exists exactly one vertex $v$ corresponding to the conjugacy class $K$. Absolute values of labels at $v$ are never $1$ because $K$ is big. Then $i(K)$ is the GCD of all labels at $v$.

Both $b_1(H), m(H), \sigma(H)$ only depend on the conjugacy class of $H$. For any special factor $H <G$ let $[H]$ be the conjugacy class of $H$ in $G$.

\begin{defi}
 We define the following complexity, which is a triple of non-negative integers, for any non-elementary GBS group $H$:
 \[
  \mathcal{C}(H)=(b_1(H),m(H), \sigma(H))
 \]
  It does not depend on the reduced graph of groups chosen to compute it nor on $\A$. We order complexities with lexicographic order.
  
  For elementary GBS groups, we define the complexity to be $(0,0,0)$.
\end{defi}

\begin{prop}\label{prop:complexite}
 Let $G$ be a GBS group and $H$ a special factor of $G$. Then 
 $\mathcal{C}(H) < \mathcal{C}(G)$.
\end{prop}
\begin{rema}
 The proposition is true for any choice of $\A$. Since a special factor with respect to $\D^\A$ is a special factor with respect to $\D^\Z$, it suffices to prove Proposition \ref{prop:complexite} for $\D^\Z$.
\end{rema}

\begin{defi}
 Let $\Gamma$ be a graph of groups and $\Gamma'$ be a subgraph of $\Gamma$. We say that $\Gamma$ is \emph{reduced with respect to $\Gamma'$} if no collapse of $\Gamma$ in the same deformation space has a subgraph with fundamental group $\pi_1(\Gamma')$.
 
 This property is a minimality condition: if $\Gamma' \subset \Gamma$, up to collapsing some edges, we can obtain $\bar \Gamma' \subset \bar \Gamma$, where $\bar \Gamma$ is reduced with respect to $\Gamma'$ and $\pi_1(\bar \Gamma') = \pi_1(\Gamma')$.
\end{defi}
\begin{rema}
 The definition implies that if an edge $e \in \Gamma$ has label $\lambda(e) = \pm 1$ and is not a loop, then $e \subset \Gamma \setminus \Gamma'$ and $o(e) \in \Gamma'$.
\end{rema}

\begin{figure}
 \centering
 \usetikzlibrary{shapes}
\begin{tikzpicture}[]

\begin{scope}[xshift=0cm, yshift=1cm]

     \draw (4,3.2) node {\textbf{At least one vertex with two edges:}};
  %  \draw[dotted] (-0.5, -10) rectangle (12,10);

    \draw[rounded corners, dashed] (0,0) rectangle (2.52,2.4);
    \draw (0.4,1.9) node {$\Gamma_H$};
    \draw (0.9,1) node {$\dots$};
    %\draw [dotted] (1.5,1) node {$\bullet$} arc (0:360: 0.3);
    \draw  (2.4,1.2) node {$\bullet$} -- (1.9, 1.7) node [right] {$p$};
    \draw [dotted] (1.4,1.8)  -- (1.9, 1.7);
    \draw [dotted] (1.5,0.5) -- (1.9, 0.7) ;
    \draw  (2.4,1.2) -- (1.9, 0.7) node [right] {$q$};
  % reste de Gamma
    \draw (2.4, 1.2) node [above right] {$1$} node {$\bullet$} -- (4,1.8) node {$\bullet$} node [above left] {$n$};
    \draw (2.4, 1.2) node [below right] {$1$} node {$\bullet$} -- (4,0.6) node {$\bullet$} node [below left] {$m$};

  %explications
    \draw [-latex] (2,-1) node [below] {not big in $\Gamma$} -- (2.4,1);
    \draw [-latex] (4.5,-1) node [below] {big in $\Gamma$} -- (4,1.7);
    \draw [-latex] (4.5,-1)  -- (4,0.5);
    \draw [-latex] (5,2) -- (6.5,2) node [text width=2cm, left=-0.5cm] {Take the subgraph};   

\end{scope}

\begin{scope}[xshift=7cm, yshift=1cm]

    \draw[rounded corners, dashed] (0,0) rectangle (2.52,2.4);
    \draw (0.4,1.9) node {$\Gamma_H$};
    \draw (0.9,1) node {$\dots$};
    %\draw [dotted] (1.5,1) node {$\bullet$} arc (0:360: 0.3);
    \draw  (2.4,1.2) node {$\bullet$} -- (1.9, 1.7) node [right] {$p$};
    \draw [dotted] (1.4,1.8)  -- (1.9, 1.7);
    \draw [dotted] (1.5,0.5) -- (1.9, 0.7) ;
    \draw  (2.4,1.2) -- (1.9, 0.7) node [right] {$q$};
  % reste de Gamma
    \draw [dotted] (2.4, 1.2) -- (4,1.8)  ;
    \draw [dotted](2.4, 1.2)  -- (4,0.6)  ;

  %explications
    \draw [-latex] (3,-0.7) node [below, text width=2.5cm] {big in $\Gamma_H$} -- (2.4,1);

\end{scope}

\begin{scope}[xshift=0cm, yshift=-6cm]
     \draw (4,3.6) node {\textbf{Only one edge on every vertex:}};

    \draw[rounded corners, dashed] (0,0) rectangle (2.52,2.4);
    \draw (0.4,1.9) node {$\Gamma_H$};
    \draw (0.9,1) node {$\dots$};
    %\draw [dotted] (1.5,1) node {$\bullet$} arc (0:360: 0.3);
    \draw  (2.4,1.2) node {$\bullet$} -- (1.9, 1.7) node [right] {$p$};
    \draw [dotted] (1.4,1.8)  -- (1.9, 1.7);
    \draw [dotted] (1.5,0.5) -- (1.9, 0.7) ;
    \draw  (2.4,1.2) -- (1.9, 0.7) node [right] {$q$};
  % reste de Gamma
    \draw (2.4, 1.2) node [below right] {$H \cap K$} node [above right] {$1$} -- (4,1.2) node {$\bullet$} node [above left] {$n$} node [below] {$K$} node [midway, above] {$e$};

        \draw [-latex] (2.1,2.6) node [above] {not big in $\Gamma$} -- (2.4,1.4);
    \draw [-latex] (4.5,2.6) node [above] {big in $\Gamma$} -- (4,1.4);

    \draw [-latex] (1, -0.2)-- (1,-1.3) node [above right] {Collapse the edge $e$};
    \draw [-latex] (5,1) -- (6.5,1) node [text width=2cm, left=-0.5cm] {Take the subgraph};

\end{scope}
\begin{scope}[xshift=0cm, yshift=-10cm]
    \draw[rounded corners, dashed] (0,0) rectangle (2.52,2.4);
    \draw (0.4,1.9) node {$\Gamma_H$};
    \draw (0.9,1) node {$\dots$};
    %\draw [dotted] (1.5,1) node {$\bullet$} arc (0:360: 0.3);
    \draw  (2.4,1.2) node {$\bullet$}  node [above right] {$K$} -- (1.9, 1.7) node [right] {$np$};
    \draw [dotted] (1.4,1.8)  -- (1.9, 1.7);
    \draw [dotted] (1.5,0.5) -- (1.9, 0.7) ;
    \draw  (2.4,1.2) -- (1.9, 0.7) node [right] {$nq$} ;
  % reste de Gamma
  %  \draw (2.4, 1.2) node [above right] {$1$} node {$\bullet$} -- (4,1.2) node {$\bullet$} node [above left] {$n$};

        \draw [-latex] (3,-0.2) node [text width=2cm , below] {becomes big} |- (2.6,1.2);
    \draw (5.5,1) node {$i(K) = gcd (np,nq)$};

    \draw [-latex] (5.5,-0.5) |- (6.7,-1);
    \draw (9,-1) node {$i(K) = n \times i(K \cap H)$};
\end{scope}

\begin{scope}[xshift=7cm, yshift=-6.5cm]
    \draw[rounded corners, dashed] (0,0) rectangle (2.52,2.4);
    \draw (0.4,1.9) node {$\Gamma_H$};
    \draw (0.9,1) node {$\dots$};
    %\draw [dotted] (1.5,1) node {$\bullet$} arc (0:360: 0.3);
    \draw  (2.4,1.2) node {$\bullet$} node [	below  right] {$K \cap H$} -- (1.9, 1.7) node [right] {$p$};
    \draw [dotted] (1.4,1.8)  -- (1.9, 1.7);
    \draw [dotted] (1.5,0.5) -- (1.9, 0.7) ;
    \draw  (2.4,1.2) -- (1.9, 0.7) node [right] {$q$};
  % reste de Gamma
    \draw [dotted] (2.4, 1.2)  -- (4,1.2) ;

        \draw [-latex] (3,2.6) node [above, text width=2.5cm] {big in $\Gamma_H$} -- (2.5,1.4);

    \draw  (2.6, -0.6) node [] {$i(K \cap H) = gcd(p,q)$};

\end{scope}

\end{tikzpicture}
 \caption{Cases when all edges in $\Gamma \setminus \Gamma_H$ are separating.} \label{fig:detailSousgraphe}
\end{figure}

\begin{proof}[Proof of proposition \ref{prop:complexite}] 
 There exists a graph of groups $\Gamma$ for $G$ such that some subgraph $\Gamma_H \subset \Gamma$ has fundamental group $H$. 
 
 We may suppose that $\Gamma$ is reduced with respect to $\Gamma_H$. In particular, $\Gamma_H$ is reduced.
 We also suppose $H$ is not an elementary GBS subgroup.
 
 If at least one of the edges in $\Gamma \setminus \Gamma_H$ is non-separating then the first Betti number of $\Gamma_H$ is strictly smaller than $b_1(\Gamma)$ so $\C(H) < \C(G)$.
 
 \medskip
 
 If all edges in $\Gamma \setminus \Gamma_H$ are separating, then each connected component of $\Gamma \setminus \Gamma_H$ is a tree attached to $\Gamma_H$ by a single vertex. Figure \ref{fig:detailSousgraphe} illustrates this case.

 In that case we first prove that the number of big vertex stabilizers (with respect to $\Amin$) cannot increase.
 
 Suppose $v \in \Gamma_H$ is a vertex whose stabilizer in $H$ is big. If no label in $\Gamma \setminus \Gamma_H$ at $v$ is $\pm 1$, then $G_v< G$ is a big vertex stabilizer in $\Gamma$.
 
 Otherwise, let $e$ be an edge with $o(e) = v$ such that $\lambda(e) = \pm 1$. The edge $e$ is not in $\Gamma_H$ since $v$ is big in $\Gamma_H$, thus $e$ is separating. Furthermore no label at $t(e)$ is $\pm 1$, so $G_{t(e)}$ is a big vertex group containing $G_v$. Thus every big stabilizer in $H$ is a subgroup of a big stabilizer in $G$.
 
 Moreover the subgroup $G_{t(e)} < G$ can contain at most one big vertex stabilizer of $H$. This implies $m(H) \leq m(G)$.
 
 \medskip
 
 Suppose at least one of the edges in $\Gamma \setminus \Gamma_H$ has both its label distinct from $\pm 1$. One of its vertices $v$ does not belong to $\Gamma'$. The stabilizer of $v$ is big (see remark \ref{rem:grosSommets}): no edge at $v$ is a loop, and all labels at $v$ are different from $\pm 1$ because $\Gamma$ is reduced with respect to $\Gamma'$. Thus $m(H) < m(G)$ so $\C(H) < \C(G)$ again.
 
 Assume that all edges in $\Gamma \setminus \Gamma_H$ are separating and have one label equal to $\pm 1$. In that case, all edges in $\Gamma \setminus \Gamma_H$ have one vertex in $\Gamma_H$, which carries the $\pm 1$ label, and a valence 1 vertex. If at least two such edges are attached to the same vertex in $\Gamma_H$ then 
 $\Gamma_H$ has strictly fewer big vertex stabilizer classes than $\Gamma$ (see figure \ref{fig:detailSousgraphe}).
 Therefore we have to deal with the case where at most one edge of $\Gamma \setminus \Gamma_H$ is attached to each vertex of $\Gamma_H$.
 
 In that case, compute $\sigma$ in $\Gamma$ and $\Gamma_H$. Only the latter is reduced. To reduce $\Gamma$ we collapse all edges in $\Gamma \setminus \Gamma_H$. Thus we get a reduced graph $\Gamma'$ which is similar to $\Gamma_H$ but at some vertices, all labels are multiplied by a factor. At such vertices the GCD of all labels is also multiplied by the factor, so the corresponding $i_{\Gamma}(K)$ is greater than $i_{\Gamma_H} (K \cap H)$. Therefore $\sigma(H) < \sigma(G)$ so $\C(H) < \C(G)$.
\end{proof}

\begin{rema}
 If $H$ is a special factor of $G$ with respect to $\D_\text{red}(G)$, then $(b_1(H), m(H)) < (b_1(G), m(G))$. However if $F$ is another special factor of $G$ for $\D_\text{red}(G)$, such that $F<H$, we do not necessarily have $(b_1(F), m(F)) < (b_1(H), m(H))$. In fact $F$ may not be a special factor of $H$ for $\D_\text{red}(H)$ (see remark \ref{rem:facteursReduits}). Therefore this simplified complexity is not helpful to study a decreasing sequence of special factors, even with respect to $\D^{\Amin}$.
\end{rema}

We deduce:
\begin{cor}\label{coro:suiteDecroissante}
 Every decreasing sequence of special factors with respect to $\D^\A$ is stationary.
\end{cor}
\begin{proof}
 Let $G \supset G_1 \supset \dots \supset G_n \supset \dots$ be a decreasing sequence of special factors of $G$.
 Either all the $G_i$ are all non-elementary groups so Proposition \ref{prop:complexite} ensures that the sequence is stationary, or for some $n \in \N$ the group $G_n$ is elementary. The elementary GBS groups do not have any non proper special factors so in that case the sequence is also stationary. 
\end{proof}

\begin{cor}\label{coro:suiteDecroissanteSysteme}
 Every $\preceq$-decreasing sequence of systems of special factors with respect to $\D^\A$ is stationary.
\end{cor}
\begin{proof}
 This is a consequence of Corollary \ref{coro:suiteDecroissante} and K\H{o}nig's Lemma. 
\end{proof}

\begin{lem}\label{lem:interDeFs}
 If $A,B$ are distinct special factors of $G$ with respect to $\D^{\mathcal A}$, then either $A\cap B$ is elliptic or it is also a special factor of $G$ with respect to $\D^{\mathcal A}$.
\end{lem}
\begin{proof}
 Let $A,B$ be as in the lemma. We will construct a tree in $\bar \D^{\mathcal A}$ in which $A \cap B$ is a vertex stabilizer.
 
 There exists a $G$-tree $T \in \bar \D^{\mathcal A}$ where $B$ is a vertex stabilizer. 
 The group $A$ acts on $T$; let $T_A$ be the minimal subtree of $T$ for this action. Subgroups of the form $hBh^{-1}\cap A$ fix a vertex in $T_A$.
 Consider the graph of groups $T_A/A$.  
 
 Let $S \in \bar \D^{\mathcal A}$ be such that $A$ is the stabilizer of a vertex $v \in S$. We may suppose that all vertices with non cyclic stabilizers are in the orbit of $v$. We perform an expansion on $S$ by replacing the vertices in the orbit of $v$ by copies of $T_A$, which is possible since edge stabilizers of $S$ are elliptic in $T_A$. We obtain a new tree $R$ which has a vertex with stabilizer $B \cap A$.
 
 We need to show that it can be obtained by collapse of a tree of $\D^{\mathcal A}$. We can perform an expansion of the vertex orbit of $S$ by replacing it by $\hat T_A$, the minimal subtree of $A$ in $\hat T \in \D^{\mathcal A}$ where $\hat T \to T$ is a collapse. All edge and vertex stabilizers of $\hat T_A$ are allowed in $\D^{\mathcal A}$. Then we get a tree $\hat R \in \D^{\mathcal A}$ which yields $R$ by collapse. This shows that $A \cap B$ is either elliptic or a special factor of $G$ with respect to $\D^{\mathcal A}$.
\end{proof}

\begin{cor}
 Let $\mathcal H, \mathcal H'$ be two systems of proper special factors. Define $\mathcal H \wedge \mathcal H':=\{[H \cap H']/ H \cap H' \text{ non elliptic }, [H] \in \mathcal H, [H'] \in \mathcal H'\}$. Then $\mathcal H \wedge \mathcal H'$ is a system of proper special factors.
\end{cor}

\begin{proof}
 The elements of $\mathcal H \wedge \mathcal H'$ are all special factors according to Lemma \ref{lem:interDeFs}. We need to check that they are all simultaneously vertex stabilizers in some tree $S \in \bar\D^\A$. The proof works like the proof of Lemma \ref{lem:interDeFs}. Let $T \in \bar \D^\A$ (resp. $T'$) be a tree in which every factor of $\mathcal H$ (resp. $\mathcal H'$) is a vertex stabilizer. We find minimal trees $T_i \subset T$ for $H_i'$ for every $[H_i'] \in \mathcal H'$, then we blow up $T'$ by replacing the vertex fixed by $gH_i'g^{-1}$ by $gT_i$ for every $g \in G$. The result is a tree $S$ in which the conjugacy classes of non-cyclic vertex stabilizers are the set $\mathcal H \wedge \mathcal H'$, and $S\in \bar \D^\A$.
\end{proof}

\begin{cor}\label{coro:facteurMing}
 Let $\G$ be a finite collection of loxodromic elements of $G$.
 The set of systems of special factors $\mathcal H$ with respect to $\D^{\mathcal A}$ such that $\G \preceq \mathcal H$ admits a smallest element for $\preceq$.
\end{cor}

\begin{proof}
 Note that this set is never empty since $\{G\}$ is itself a system of special factors with respect to $\D^{\mathcal A}$.
 
 Corollary \ref{coro:suiteDecroissanteSysteme} ensures that any $\preceq$-decreasing sequence of systems of special factors is stationary. Thus there exists a system of special factors $\mathcal H$ such that $\G \preceq \mathcal H$ which is minimal for this property.

 Let us show that it is unique. Let $\mathcal H'$ be another minimal special factor such that $\G \preceq \mathcal H'$. Let $g \in \G$: there is $H \in \mathcal H$ and $H' \in \mathcal H'$ such that $g \in H \cap H'$ so $\mathcal H \wedge \mathcal H'$ is not empty. By Lemma \ref{lem:interDeFs} it is a system of special factors of $G$.
 
 By minimality of $\mathcal H$ we get $\mathcal H \subset \mathcal H'$, and conversely by minimality of $\mathcal H'$. By Remark \ref{rem:ordre} 2. the relation $\preceq$ is an order so both systems are equal.
\end{proof}

Using Whitehead algorithm and corollary \ref{coro:facteurMing} we deduce that there exists an algorithm which finds the smallest system of special factors $\mathcal H$ such that $\G \subset G$ is $\mathcal H$-peripheral with respect to $\D^{\mathcal A}$:
\begin{theo}
 There is an algorithm which takes as input
 \begin{itemize}
  \item a marked graph of groups $\Gamma_G$ representing a group $G$
  \item a finite collection of loxodromic elements $\G \subset G$ as loops in $\Gamma_G$
  \item finite sets $(I_v)_{v \in V(\Gamma_G)}$ representing a collection of allowed edge groups $\mathcal A$
 \end{itemize}
 and outputs the smallest system of special factors $\mathcal H$ of $G$ with respect to $\D^{\mathcal A}$ such that $\G$ is $\mathcal H$-peripheral, as a finite collection of marked graph of groups.
\end{theo}
\begin{proof}
 We give the algorithm in the case where $\G$ has a single element. The algorithm consists in constructing a decreasing sequence of special factors containing $g$. In the case where $\G$ consists of more than one element, the algorithm would construct a decreasing sequence of systems of special factors.

 Define $\mathcal H_0:=\{G\}$, $\Gamma_0 := \Gamma_G$ and $\D_0:=\D^{\mathcal A}$. 
 We will construct a decreasing sequence of special factors $H_0 \supset \dots \supset H_N$, where $g \in H_n$ for every $n \in \{0, \dots, N\}$. The space $\D_n$ is defined as the induced deformation space for $H_n$. It coincides with the induced deformation space for $H_n$ seen as a special factor of some $H_{m}$ for $m<n$.
 For every $n$, we will construct sets $(I^n_v)_{v \in V(\Gamma_n)}$ which represent the family of allowed edge groups of $\D_n$.
 
 Here is the algorithm. Start with $i=0$. Use Theorem \ref{theo:algoWhitehead} applied to $\Gamma_i$ to decide whether $g$ is contained in a special factor of $H_i$ or not. If yes, then $H_i$ is the minimal special factor containing $g$ and the algorithm stops. Else the algorithm gives 
 \begin{itemize}
   \item a special factor $H_{i+1} \subsetneq H_i$ such that $g \in H_{i+1}$
   \item a new graph of groups $\Gamma'_i$ with $\pi_1(\Gamma'_i) \simeq H_i$ 
   \item a subgraph of groups $\Gamma_{i+1} \subset \Gamma'_i$ such that $\pi_1(\Gamma_{i+1}) \simeq H_{i+1}$.
 \end{itemize}
 Then start again with $i+1$ instead of $i$.
 
 Applying this construction, we get a decreasing sequence of special factors with respect to $\D^{\A}$.

Corollary \ref{coro:suiteDecroissante} guarantees that this sequence is stationary, which means that the algorithm stops eventually.  
The last factor obtained, given as a marked graph of groups, is the minimal special factor containing $g$.
\end{proof}

\bibliography{bibli}

\end{document}